\documentclass[leqno,11pt]{amsart}
\usepackage{amsmath,amscd,amsthm,amsxtra}
\usepackage{epsfig,graphics,color,colortbl}
\usepackage{amssymb,latexsym}
\usepackage{mathrsfs}
\usepackage{eucal,upgreek}
\usepackage[poly,all]{xy}
\usepackage{hyperref}
\usepackage{tikz-cd}

\newtheorem{thm}{\bf Theorem}[section]

\newtheorem{prop}[thm]{\bf Proposition}
\newtheorem{cor}[thm]{\bf Corollary}

\newtheorem{rem}[thm]{\bf Remark}
\newtheorem{ex}[thm]{\bf Example}

\numberwithin{equation}{section}

\newcommand{\bs}{\boldsymbol}
\newcommand{\A}{\mathbb{A}}

\newcommand{\bB}{\mathbb{B}}

\newcommand{\cB}{\mathcal{B}}
\newcommand{\W}{\mathcal{W}}
\newcommand{\cP}{\mathscr{P}}

\newcommand{\pf}{\noindent{\bfseries Proof. }}
\newcommand{\ov}{\overline}

\newcommand{\ba}{\bs{\rm a}}
\newcommand{\bb}{\bs{\rm b}}

\newcommand{\bi}{\bs{\rm i}}

\newcommand{\M}{{\mathcal{M}}}

\newcommand{\cM}{\mathcal{M}}

\newcommand{\N}{\mathbb{N}}
\newcommand{\gl}{\mathfrak{gl}}
\newcommand{\Z}{\mathbb{Z}}
\newcommand{\C}{\mathbb{C}}

\newcommand{\te}{\widetilde{e}}
\newcommand{\tf}{\widetilde{f}}
\newcommand{\g}{\mathfrak{g}}

\newcommand{\mc}{\mathcal}
\newcommand{\mf}{\mathfrak}

\newcommand{\la}{\lambda}

\newcommand{\nw}{^{\nwarrow}}

\begin{document}
\title[Lusztig data of Young tableaux]
{A crystal embedding into Lusztig data of type $A$}
\author{JAE-HOON KWON}

\address{Department of Mathematical Sciences, Seoul National University, Seoul 08826, Korea}
\email{jaehoonkw@snu.ac.kr}


\thanks{This work was supported by Samsung Science and Technology Foundation under Project Number SSTF-BA1501-01.}

\begin{abstract}
Let $\bi$ be a reduced expression of the longest element in the Weyl group of type $A$, which is adapted to a Dynkin quiver with a single sink. 
We present a simple description of the crystal embedding of Young tableaux of arbitrary shape into $\bi$-Lusztig data, which also gives an algorithm for the transition matrix between Lusztig data associated to reduced expressions adapted to quivers with a single sink. 
\end{abstract}

\maketitle

\section{Introduction}

Let $U_q(\g)$ be the quantized enveloping algebra associated to a symmetrizable Kac-Moody algebra $\g$. The negative part of $U_q(\g)$ has a basis called a canonical basis \cite{Lu91} or lower global crystal basis \cite{Kas91}, which has many fundamental properties. 
The canonical basis forms a colored oriented graph $B(\infty)$, called a crystal, with respect to Kashiwara operators. The crystal $B(\infty)$ plays an important role in the study of combinatorial aspects of $U_q(\g)$-modules together with its subgraph $B(\la)$ associated to any integrable highest weight module $V(\la)$ with highest weight $\la$.

Suppose that $\g$ is a finite-dimensional semisimple Lie algebra with the index set $I$ of simple roots. Let $\bi=(i_1,\ldots,i_N)$ be a sequence of indices in $I$ corresponding to a reduced expression of the longest element in the Weyl group of $\g$. 
A PBW basis associated to $\bi$ is a basis of the negative part of $U_q(\g)$ \cite{Lu93}, which is parametrized by the set $\cB_{\bi}$ of $N$-tuple of non-negative integers. 
One can identify $B(\infty)$ with $\cB_{\bi}$ since the associated PBW basis coincides with the canonical basis at $q=0$ \cite{S94}. We call an element in $\cB_{\bi}$ an $\bi$-Lusztig datum or Lusztig parametrization associated to $\bi$.

Consider the map
\begin{equation}\label{eq:embedding_psi}
\begin{split}
\xymatrixcolsep{3pc}\xymatrixrowsep{4pc}\xymatrix{
\psi_\la^{\bi} :\ B(\la)\otimes T_{-\la}  \ \ar@{^{(}->}[r]  &\ \cB_{\bi}},
\end{split}
\end{equation}
given by the $\bi$-Lusztig datum of $b\in B(\la)$ under the embedding of $B(\la)\otimes T_{-\la}$ into $B(\infty)$, where  $T_{-\la}=\{t_{-\la}\}$ is an abstract crystal with ${\rm wt}(t_{-\la})=-\la$ and $\varphi_i(t_{-\la})=-\infty$ for $i\in I$. 
In this paper, we give a simple combinatorial description of \eqref{eq:embedding_psi} when $\g=\mf{gl}_n$ and $\bi$ is a reduced expression adapted to a Dynkin quiver of type $A_{n-1}$ with a single sink (Theorem \ref{thm:main-1}). It is well-known that when $\bi$ is adapted to a quiver with one direction, for example $\bi=\bi_0=(1,2,1,3,2,1,\ldots,n-1,\ldots,1)$, the $\bi_0$-Lusztig datum of a Young tableaux is simply given by counting the number of occurrences of each entry in each row. But the $\bi$-Lusztig datum for arbitrary $\bi$ is not easy to describe in general, and one may apply a sequence of Lusztig's transformations \cite{Lu93} or the formula for a transition map $R_{\bi_0}^{\bi} : \cB_{\bi_0}\rightarrow \cB_{\bi}$ by Berenstein-Fomin-Zelevinsky \cite{BFZ}. 
  We remark that our algorithm for computing $\psi_\la^{\bi}$ is completely different from the known methods, and hence provides an alternative description of $R^{\bi}_{\bi_0}$. 

Let us explain the basic ideas in our description of $\psi_\la^{\bi}$.
Suppose that $\Omega$ is a quiver of type $A_{n-1}$ with a single sink and $\bi$ is adapted to $\Omega$.
Let $J\subset I$ be a maximal subset such that each connected component of the corresponding quiver $\Omega_J\subset \Omega$ has only one direction. Let $\g_J$ be the maximal Levi subalgebra and $\mf{u}_J$ the nilradical associated to $J$, respectively.

The first step is to prove a tensor product decomposition $\cB_{\bi}\cong B^J(\infty)\otimes B_J(\infty),$
as a crystal,
where $B_J(\infty)$ is the crystal of the negative part of $U_q(\g_J)$ and $B^J(\infty)$ is the crystal of the quantum nilpotent subalgebra $U_q(\mf{u}_J)$. The isomorphism is just given by restricting the Lusztig datum to each part, and it is a special case of the bijection introduced in \cite{BKT,S94} using crystal reflections. Here we show that it is indeed a morphism of crystals by using Reineke's description of $B(\infty)$ in terms of representations of $\Omega$ \cite{Re}. We refer the reader to a recent work by Salibury-Schultze-Tingley \cite{SST} on $\bi$-Lusztig data, {which also implies the combinatorial description of Kashiwara operators on $B(\infty)$ used in this paper.}

The next step is to construct an embedding of $B(\la)\otimes T_{-\la}$ into $B^J(\infty)\otimes B_J(\infty)$ using a crystal theoretic interpretation of Sagan and Stanley's skew RSK algorithm \cite{SS}, which was observed in the author's previous work \cite{K09} (see also \cite{K13,K14}), and using the embedding \eqref{eq:embedding_psi} in case of $\bi$ adapted to a quiver with one direction. Hence we obtain an $\bi$-Lusztig
datum of a Young tableau for any $\bi$ adapted to $\Omega$. {One may consider the image of the embedding by using a combinatorial description of $\ast$-crystal structure on $\cB_{\bi}$ in \cite{Re}, but we do not discuss it here.}

{Our description of the embedding $\psi_\la^{\bi}$ also provides an algorithm for a transition map $R_{\bi_0}^{\bi} : \cB_{\bi_0}\rightarrow \cB_{\bi}$ together with its inverse $R_{\bi}^{\bi_0}$ since $\psi_\la^{\bi}$ naturally extends to an isomorphism from another realization of $B(\infty)$ given by the set of {large tableaux} \cite{CL,HL}. Therefore we obtain an algorithm for a transition map $R^{\bi'}_{\bi}=R^{\bi'}_{\bi_0}\circ R^{\bi_0}_{\bi}$ for any $\bi$ and $\bi'$ which are adapted to quivers with a single sink.} Roughly speaking, $R^{\bi'}_{\bi}$ is given by a composition of skew RSK and its inverse algorithms with respect to various maximal Levi subalgebras depending on $\bi$ and $\bi'$.


The paper is organized as follows: In Sections \ref{sec:prel} and \ref{sec:yt}, we review necessary background on crystals and related materials. In Section \ref{sec:pbw}, we give an explicit description of the crystal $\cB_{\bi}$ when $\bi$ is adapted to a Dynkin quiver of type $A_{n-1}$ with a single sink, and then prove the decomposition of $B_{\bi}$ as a tensor product of two subcrystals. Finally in Section \ref{sec:embedding}, we construct an embedding of the crystal of Young tableaux of arbitrary shape $\la$ into $\cB_{\bi}$.\vskip 2mm

{\bf Acknowledgement} The author would like to thank Myungho Kim for valuable discussion and kind explanation on representations of quivers. 

\section{Review on crystals}\label{sec:prel}

\subsection{}
Let us give a brief review on crystals (see \cite{HK,Kas91,Kas95} for more details). We denote by $\Z_+$ the set of non-negative integers.
Fix a positive integer $n$ greater than $1$. Throughout the paper, $\g$ denotes the general linear Lie algebra $\gl_n(\C)$ which is spanned by the elementary matrices $e_{ij}$ for $1\leq i,j\leq n$. Let $P^\vee=\bigoplus_{i=1}^n\Z e_{ii}$ be the dual weight lattice and $P={\rm Hom}_\Z(P^\vee,\Z)=\bigoplus_{i=1}^n\Z\epsilon_i$ be the weight lattice of $\g$ with $\langle \epsilon_i,e_{jj} \rangle =\delta_{ij}$ for $i,j$. Define a symmetric bilinear form $(\, \cdot \,|\, \cdot\,)$ on $P$ such that $(\epsilon_i |\epsilon_j)=\delta_{ij}$ for $i,j$. 
Set $I=\{1,\ldots, n-1\}$. Then $\{\,\alpha_i:=\epsilon_i-\epsilon_{i+1}\,|\,i\in I\,\}$ is the set of simple roots and $\{\,h_i:=e_{ii}-e_{i+1\,i+1}\,|\,i\in I\,\}$ is the set of simple coroots of $\g$. Let $\Phi^+=\{\,\epsilon_i-\epsilon_j\,|\,1\leq i<j\leq n\,\}$ denote the set of positive roots of $\g$.

Let $W\cong\mf S_n$ be the Weyl group of $\g$, which is generated by simple reflections $s_i$ for $i\in I$. Let $w_0$ be the longest element in $W$, which is of length  $N:=n(n-1)/2$,
and let $R(w_0)=\{\,(i_1,\ldots,i_N)\,|\,w_0=s_{i_1}\ldots s_{i_N}\,\}$ be the set of reduced expressions of $w_0$.  

For $J\subset I$, let $\g_J$ be the subalgebra of $\g$ generated by $e_{ii}$ for $1\leq i\leq n$ and the root vectors associated to $\pm \alpha_j$ for $j\in J$. Let $\Phi^+_J$ be the set of positive roots of $\g_J$ and $\Phi^+(J)=\Phi^+\setminus \Phi^+_J$.

A {\it $\g$-crystal} is a set
$B$ together with the maps ${\rm wt} : B \rightarrow P$,
$\varepsilon_i, \varphi_i: B \rightarrow \mathbb{Z}\cup\{-\infty\}$ and
$\te_i, \tf_i: B \rightarrow B\cup\{{\bf 0}\}$ for $i\in I$ satisfying the following conditions:  for $b\in B$ and $i\in I$,
\begin{itemize}
\item[(1)]  
$\varphi_i(b) =\langle {\rm wt}(b),h_i \rangle +
\varepsilon_i(b)$,

\item[(2)] $\varepsilon_i(\te_i b) = \varepsilon_i(b) - 1,\ \varphi_i(\te_i b) =
\varphi_i(b) + 1,\ {\rm wt}(\te_ib)={\rm wt}(b)+\alpha_i$ if $\te_i b \in B$,

\item[(3)] $\varepsilon_i(\tf_i b) = \varepsilon_i(b) + 1,\ \varphi_i(\tf_i b) =
\varphi_i(b) - 1,\ {\rm wt}({\tf_i}b)={\rm wt}(b)-\alpha_i$ if $\tf_i b \in B$,

\item[(4)] $\tf_i b = b'$ if and only if $b = \te_i b'$ for $b' \in B$,

\item[(5)] $\te_ib=\tf_ib={\bf 0}$ when $\varphi_i(b)=-\infty$.
\end{itemize}
Here ${\bf 0}$ is a formal symbol and $-\infty$ is the smallest
element in $\Z\cup\{-\infty\}$ such that $-\infty+n=-\infty$
for all $n\in\Z$. Unless otherwise specified, a crystal means a $\g$-crystal throughout the paper for simplicity. 

Let $B_1$ and $B_2$ be crystals.
A {\em tensor product $B_1\otimes B_2$}
is a crystal, which is defined to be $B_1\times B_2$  as a set with elements  denoted by
$b_1\otimes b_2$, where  
{\allowdisplaybreaks
\begin{equation*}
\begin{split}
{\rm wt}(b_1\otimes b_2)&={\rm wt}(b_1)+{\rm wt}(b_2), \\
\varepsilon_i(b_1\otimes b_2)&= {\rm
max}\{\varepsilon_i(b_1),\varepsilon_i(b_2)-\langle {\rm
wt}(b_1),h_i\rangle\}, \\
\varphi_i(b_1\otimes b_2)&= {\rm max}\{\varphi_i(b_1)+\langle {\rm
wt}(b_2),h_i\rangle,\varphi_i(b_2)\},
\end{split}
\end{equation*}}
{\allowdisplaybreaks
\begin{equation}\label{eq:tensor product of crystals}
\begin{split}
{\te}_i(b_1\otimes b_2)&=
\begin{cases}
{\te}_i b_1 \otimes b_2, & \text{if $\varphi_i(b_1)\geq \varepsilon_i(b_2)$}, \\
b_1\otimes {\te}_i b_2, & \text{if
$\varphi_i(b_1)<\varepsilon_i(b_2)$},
\end{cases}\\
{\tf}_i(b_1\otimes b_2)&=
\begin{cases}
{\tf}_i b_1 \otimes b_2, & \text{if  $\varphi_i(b_1)>\varepsilon_i(b_2)$}, \\
b_1\otimes {\tf}_i b_2, & \text{if $\varphi_i(b_1)\leq
\varepsilon_i(b_2)$},
\end{cases}
\end{split}
\end{equation}}
\noindent for $i\in I$. Here we assume that ${\bf 0}\otimes
b_2=b_1\otimes {\bf 0}={\bf 0}$.

A {\it morphism}
$\psi : B_1 \rightarrow B_2$ is a map from $B_1\cup\{{\bf 0}\}$ to
$B_2\cup\{{\bf 0}\}$ such that
\begin{itemize}
\item[(1)] $\psi(\bf{0})=\bf{0}$,

\item[(2)] ${\rm wt}(\psi(b))={\rm wt}(b)$,
$\varepsilon_i(\psi(b))=\varepsilon_i(b)$, and
$\varphi_i(\psi(b))=\varphi_i(b)$ when $\psi(b)\neq \bf{0}$,

\item[(3)] $\psi(\te_i b)=\te_i\psi(b)$ when $\psi(b)\neq \bf{0}$ and
$\psi(\te_i b)\neq \bf{0}$,

\item[(4)] $\psi(\tf_i
b)=\tf_i\psi(b)$ when $\psi(b)\neq \bf{0}$ and $\psi(\tf_i b)\neq \bf{0}$,
\end{itemize}
for $b\in B_1$ and $i\in I$.
We call $\psi$ an {\it embedding} and $B_1$ a {\it subcrystal of}
$B_2$ when $\psi$ is injective.

The {\it dual crystal $B^\vee$} of a crystal $B$ is defined
to be the set $\{\,b^\vee\,|\,b\in B\,\}$ with ${\rm
wt}(b^\vee)=-{\rm wt}(b)$, $\varepsilon_i(b^\vee)=\varphi_i(b)$,
$\varphi_i(b^\vee)=\varepsilon_i(b)$, $\te_i(b^\vee)=(\tf_i b)^\vee$, and $\tf_i(b^\vee)=\left(\te_i b \right)^\vee$ for
$b\in B$ and $i\in I$. We assume that ${\bf 0}^\vee={\bf 0}$.

For $\mu\in P$, let $T_\mu=\{t_\mu\}$ be a crystal, where ${\rm wt}(t_\mu)=\mu$, $\te_it_\mu=\tf_it_\mu={\bf 0}$, and $\varepsilon_i(t_\mu)=\varphi_i(t_\mu)=-\infty$ for all $i\in I$.

\subsection{}

Let $q$ be an indeterminate. Let $U=U_q(\g)$ be the quantized enveloping algebra of $\g$, which is an associative $\mathbb{Q}(q)$-algebra with $1$ generated by $e_i$, $f_i$, and $q^h$ for $i\in I$ and $h\in P^\vee$. 
Let $U^-=U^-_q(\g)$ be the negative part of $U$, the subalgebra generated by $f_i$ for $i\in I$.
We put  $[m]=\frac{q^m-q^{-m}}{q-q^{-1}}$ and $[m]!=[1][2]\cdots [m]$ for $m\in\N$. Let $t_i=q^{h_i}$, $e_i^{(m)}=e_i^m/[m] !$, and $f_i^{(m)}=f_i^m/[m]!$  for $m\in\N$ and $i\in I$. 
Let $A_0$ denote the subring of $\mathbb{Q}(q)$ consisting
of rational functions regular at $q=0$.

For $i\in I$, let $T_i$ be the $\mathbb{Q}(q)$-algebra automorphism of $U$ given by
{\allowdisplaybreaks
\begin{align*}
T_i(t_j)&=t_j t_i^{-a_{ij}}, \\
T_i(e_j)&=  
\begin{cases}
-f_it_i, & \text{if $j=i$}, \\
\sum_{k+l=-a_{ij}} (-1)^kq^{-l}e_i^{(k)}e_je_i^{(l)}, & \text{if $j\neq i$},
\end{cases}\\
T_i(f_j)&=  
\begin{cases}
-t_i^{-1}e_i, & \text{if $j=i$}, \\
\sum_{k+l=-a_{ij}} (-1)^kq^{l}f_i^{(k)}f_jf_i^{(l)}, & \text{if $j\neq i$},
\end{cases}
\end{align*}}
for $j\in I$, where $a_{ij}=\langle \alpha_j,h_i \rangle$.
Note that $T_i$ is denoted by {$T''_{i,1}$} in \cite{Lu93} (see also \cite{S94}). 

For ${\bi}=(i_1,\ldots, i_N)\in R(w_0)$ and ${\bf c}=(c_1,\ldots,c_N)\in\Z_+^N$, consider the vectors of the following form:
\begin{equation}\label{eq:PBW vector}
\begin{split}
b _{\bi}(\bf c)=&
f^{(c_{1})}_{i_1}T_{i_1}(f^{(c_{2})}_{i_2})\cdots T_{i_{1}}T_{i_2}\cdots T_{i_{N-1}}(f^{(c_{N})}_{i_{N}}).
\end{split}
\end{equation}
The set $B_{\bf i}:=\{\,b _{\bi}({\bf c})\,|\,{\bf c}\in\Z_+^{N}\,\}$ is a $\mathbb{Q}(q)$-basis of $U^-$, which is often referred to as a {\it PBW basis} \cite{Lu93}.

The $A_0$-lattice of $U^-$ generated by $B_{\bf i}$ is independent of the choice of ${\bf i}$, which we denote by $L(\infty)$. If $\pi : L(\infty) \rightarrow L(\infty)/q L(\infty)$ is the canonical projection, then $\pi(B_{\bf i})$ is a $\mathbb{Q}$-basis of  $L(\infty)/qL(\infty)$ and also independent of the choice of ${\bf i}$, which we denote by $B(\infty)$. 
Indeed the pair $(L(\infty),B(\infty))$ coincides with the {\it Kashiwara's crystal base of $U^-$} \cite{Kas91}, that is, $L(\infty)$ is invariant under $\te_i$, $\tf_i$, and  $\te_iB(\infty)\subset B(\infty)\cup \{0\}$, $\tf_iB(\infty)\subset B(\infty)\cup \{0\}$ for $i\in I$, where $\te_i$ and $\tf_i$ denote the modified Kashiwara operators on $U^-$ given by
\begin{equation*}
\te_i x =\sum_{k\geq 1}f_i^{(k-1)}x_k,\quad\quad \tf_i x =\sum_{k\geq 0}f_i^{(k+1)}x_k,
\end{equation*}
for $x=\sum_{k\geq 0}f_i^{(k)}x_k\in U^-$, where $x_k\in T_i(U^-)\cap U^-$ for $k\geq 0$ (see \cite{Lu96,S94}).
The set $B(\infty)$ equipped with the induced operators $\te_i$ and $\tf_i$ becomes a crystal, where $\varepsilon_i(b)=\max\{\,k\,|\,\te_i^kb\neq 0\,\}$ for $i\in I$ and $b\in B(\infty)$. 

Let $P^+=\{\,\la\in P\,|\, \langle \la,h_i\rangle \geq 0\ \text{for $i\in I$} \,\}$ be the set of dominant integral weights. For $\la\in P^+$, let $V(\la)$ be the irreducible highest weight $U$-module with highest weight $\la$, which is given by $U^-/\sum_{i\in I}U^- f_i^{\langle \la,h_i \rangle +1}\cdot 1$ as a left $U^-$-module. If $\pi_\la : U^- \rightarrow V(\la)$ is the canonical projection, then $L(\la):=\pi_\la(L(\infty))$ is an $A_0$-lattice of $V(\la)$ and $B(\la): =\pi_\la(B(\infty))\setminus\{0\}$ is a $\mathbb{Q}$-basis of $L(\la)/qL(\la)$.  The pair $(L(\la),B(\la))$ is called the {\em crystal base of $V(\la)$}. The set $B(\la)$ becomes a crystal with respect to $\te_i$ and $\tf_i$ induced from those on $B(\infty)$, where $\varepsilon_i(b)=\max\{\,k\,|\,\te_i^kb\neq 0\,\}$ and $\varphi_i(b)=\max\{\,k\,|\,\tf_i^kb\neq 0\,\}$ for $i\in I$ and $b\in B(\la)$ \cite{Kas91}.

\section{Crystal of Young tableaux}\label{sec:yt}

\subsection{} 

Let us recall some necessary background on semistandard tableaux and related combinatorics following \cite{Ful}. 
Let $\mathscr{P}$ be the set of partitions. 
We identify $\lambda=(\lambda_i)_{i\geq 1}\in \cP$ with a {Young diagram}. Let $\la/\mu$ denote a skew Young diagram associated to $\la,\mu\in \cP$ with $\la\supset \mu$, and let $(\lambda/\mu)^\pi$ denote the skew Young diagram obtained by $180^{\circ}$-rotation of $\lambda/\mu$.
 
Let $\A$ be a linearly ordered set. For a skew Young diagram
$\lambda/\mu$,  let $SST_\A(\lambda/\mu)$ be the set of all
semistandard tableaux of shape $\lambda/\mu$ with entries in $\A$.
 Let $\W_{\A}$ be the set of finite words in
$\A$. For $T\in SST_\A(\lambda/\mu)$, let ${\rm sh}(T)$ denote the shape of $T$, 
and let $w(T)$ be a word in $\W_\A$ obtained by reading the entries of $T$ row
by row from top to bottom, and from right to left in each row.

Let $T\in SST_{\A}(\la^\pi)$ be given for $\la\in \cP$. For $a\in \A$, we define $T \leftarrow a$ to be the tableau obtained by applying the Schensted's column insertion of $a$ into $T$ in a reverse way starting from the rightmost column of $T$ so that ${\rm sh}(T\leftarrow a)=\mu^\pi$ for some $\mu\supset \la$ obtained by adding a box in a corner of $\la$. 
We also denote by $T^{\nw}$ the unique tableau in $SST_{\A}(\la)$, which is Knuth equivalent to $T$. 
{Note that the map $T\mapsto T^{\nw}$ gives a bijection from $SST_{\A}(\la^\pi)$ to $SST_\A(\la)$, where the inverse map is given by $((\cdots(a_r\leftarrow a_{r-1})\leftarrow \cdots)\leftarrow a_1)$ for $S\in SST_{\A}(\la)$ with $w(S)=a_1\cdots a_r$. }

Let $\bB$ be another linearly ordered set, and let
\begin{equation*}\label{MAB}
\M_{\A\times \bB}=\left\{\,M=(m_{ab})_{a\in\A,b\in \bB}\,\,\Bigg\vert\,\, m_{ab}\in\Z_{+},\ \ \sum_{a,b}m_{ab}<\infty\,\right\}.
\end{equation*}
Let $\mathcal{I}_{\A\times\bB}$ be the set of biwords $(\ba,\bb)\in
\W_{\A}\times \W_{\bB}$ such that (1) $\ba=a_1\cdots a_r$ and
$\bb=b_1\cdots b_r$  for some $r\geq 0$, (2)   $(a_1,b_1)\leq \cdots
\leq (a_r,b_r)$, where for $(a,b)$ and $(c,d)\in \A\times \bB$,
$(a,b)< (c,d)$ if and only if $(a<c)$ or ($a=c$ and $b>d$). 
There is a bijection  
\begin{equation}\label{eq:bijection I to M}
\begin{split} 
&\xymatrixcolsep{3pc}\xymatrixrowsep{0pc}\xymatrix{
\mathcal{I}_{\A\times\bB} \ \ar@{->}[r]  & \ \M_{\A\times\bB}   \\
(\ba,\bb) \ \ar@{|->}[r]  & \ M(\ba,\bb) }
\end{split}
\end{equation}
where
$M(\ba,\bb)=(m_{ab})$ with
$m_{ab}=\left|\{\,k\,|\,(a_k,b_k)=(a,b) \,\}\right|$ and the pair of
empty words $(\emptyset,\emptyset)$ corresponds to the zero matrix $O$. 

For $(\ba,\bb)\in \mc{I}_{\A\times \bB}$,
we write $M[\bb,\ba]=M(\ba,\bb)^t\in \M_{\bB\times\A}$, where $M^t$ denotes the transpose of $M\in \M_{\bB\times\A}$. 
For $(\ba,\bb)\in \mc{I}_{\A\times\bB}$, there exist unique 
$\ba^\tau\in \W_{\A}$ and 
$\bb^\tau\in \W_{\bB}$, which are rearrangements of $\ba$ and $\bb$, respectively, satisfying $M(\bb^\tau,\ba^\tau) = M(\ba,\bb)^t \in \M_{\bB\times \A}$ with $(\bb^\tau,\ba^\tau)\in \mc{I}_{\bB\times \A}$, or equivalently
\begin{equation}\label{eq:Mtau}
M[\ba^\tau,\bb^\tau] = M(\ba,\bb)\in \M_{\A\times \bB}. 
\end{equation}

Fix $\la\in\cP$. Let $T\in SST_{\A}(\la^\pi)$ and $M\in \M_{\A\times\bB}$ be given, 
where 
$M=M(\ba, \bb)$ for some $(\ba,\bb)\in \mathcal{I}_{\A\times \bB}$. 
Suppose that 
$\ba^\tau=a^\tau_1\cdots a^\tau_r$ and $\bb^\tau=b^\tau_1\cdots b^\tau_r$. 
We define the pair of tableaux ${\bf P}(T \leftarrow M)$ and ${\bf Q}(T \leftarrow M)$ inductively as follows:
For $1\leq i\leq r$, 
put ${\bf P}^{(i)}= ({\bf P}^{(i-1)}\leftarrow a^\tau_{r-i+1})$, and $\la^{(i)}={\rm sh}\left({\bf P}^{(i)}\right)^\pi$ with ${\bf P}^{(0)}=T$ and $\la^{(0)}=\la$. 
Define ${\bf P}(T \leftarrow M)= {\bf P}^{(r)}$ and $\mu={\rm sh}\left({\bf P}^{(r)}\right)^\pi$, and define ${\bf Q}(T \leftarrow M)$ to be the tableau of shape $\left(\mu/\la \right)^\pi$, where $\left(\la^{(i)}/\la^{(i-1)}\right)^\pi$ is filled with $b^\tau_{r-i+1}$ for $1\leq i\leq r$.
Then the map
\begin{equation}\label{eq:skew_RSK_antinormal}
\begin{split}
\xymatrixcolsep{3pc}\xymatrixrowsep{0.5pc}\xymatrix{
\kappa : SST_{\A}(\la^\pi) \times \M_{\A\times\bB} \ \ar@{->}[r]  & \ \
\displaystyle\bigsqcup_{\mu\supset \la} SST_{\A}(\mu^\pi)\times SST_{\bB}(\left(\mu/\la\right)^\pi)  \\
 (T,M)\  \ar@{|->}[r]  &\ \ ({\bf P}(T \leftarrow M),{\bf Q}(T \leftarrow M))}
\end{split}
\end{equation}
is a bijection, which is a skew analogue of the usual RSK correspondence \cite{SS}.

\subsection{}

Let $[n]=\{\,1<\cdots<n\,\}$ and $[\ov{n}]=\{\,\ov{n}<\cdots< \ov{1}\,\}$ be linearly ordered sets.
We regard $[n]$ as a crystal $B(\epsilon_1)$, where ${\rm wt}(k)=\epsilon_k$ for $k\in[n]$, and $[\ov{n}]$ as the dual crystal $[n]^\vee$, where $\ov{k}=k^\vee$ for $k\in [n]$.
Then $\W_{[n]}$ and $\W_{[\ov{n}]}$ are crystals, where
we identify $w=w_1\ldots w_r$ with $w_1\otimes \cdots \otimes w_r$. 
The crystal structure on $\W_{[n]}$ is easily described by so-called {\em signature rule} 
(cf. \cite[Section 2.1]{KashNaka}).

Let $\cP_n$ be the set of partitions $\la=(\la_1,\ldots,\la_n)$ of length less than or equal to $n$. For $\la\in \cP_n$, $SST_{[n]}(\la)$ is a crystal under the identification of $T$ with $w(T)\in \W_{[n]}$, and it is isomorphic to $B(\la)$, where we regard $\la$ as $\sum_{i=1}^n\la_i\epsilon_i\in P^+$ \cite{KashNaka}, while $SST_{[\ov{n}]}(\la)$ is isomorphic to $B(-w_0\la)$. One can define a crystal structure on $SST_{[n]}(\mu/\nu)$ for a skew Young diagram $\mu/\nu$ in a similar way.
Note that $SST_{[n]}(\la)^\vee \cong SST_{[\ov{n}]}(\la^\pi)$, where the isomorphism is given by taking the $180^{\circ}$-rotation and replacing the entry $i$ with $\ov{n-i+1}$ for $i\in [n]$.

For $0\leq t\leq n$, let 
$\sigma^{-1} : SST_{[n]}(1^t) \longrightarrow SST_{[\ov{n}]}(1^{n-t})$
be a bijection, where $\sigma^{-1}(T)$ is the tableau with entries $[\ov{n}] \setminus \{\ov{k_1}, \ldots, \ov{k_t} \}$ for $T$ with entries $k_1< \cdots <  k_t$.
For $d\geq \la_1$, 
define 
\begin{equation}\label{eq:tau}
\xymatrixcolsep{3pc}\xymatrixrowsep{4pc}\xymatrix{
 \sigma^{-d} : SST_{[n]}(\la) \ \ar@{->}[r]  & \ SST_{[\ov{n}]}(\sigma^{-d}(\la))}, 
\end{equation}
where $\sigma^{-d}(\la)=(d^n)/\la$, and the $i$th column of $\sigma^{-d}(T)$ from the left is obtained by applying $\sigma^{-1}$ to the $i$th column of $T\in SST_{[n]}(\la)$ (which is assumed to be empty if $i>\la_1$). Then $\sigma^{-d}$ commutes with $\te_i$ and $\tf_i$ for $i\in I$, where ${\rm wt}(\sigma^{-d}(T))={\rm wt}(T)-d(\epsilon_1+\cdots+\epsilon_n)$. 
We have an isomorphism of crystals
\begin{equation}\label{eq:tau_iso}
\xymatrixcolsep{3pc}\xymatrixrowsep{4pc}\xymatrix{
SST_{[n]}(\la)\otimes T_{\xi} \ \ar@{->}[r]  & \ SST_{[\ov{n}]}(\sigma^{-d}(\la))},
\end{equation} 
where $\xi=-d(\epsilon_1+\cdots+\epsilon_n)$. Also \eqref{eq:tau} and \eqref{eq:tau_iso} hold when $[n]$ and $[\ov{n}]$ are exchanged (assuming that $\ov{\ov{k}}=k$ for $k\in [n]$).

\section{Crystal of Lusztig  data}\label{sec:pbw}

\subsection{}
Let ${\bf i}=(i_1,\ldots,i_N)\in R(w_0)$ be given. 
We have 
\begin{equation*}
\Phi^+=\{\,\beta_1:=\alpha_{i_1},\  \beta_2:=s_{i_1}(\alpha_{i_2}), \ \ldots ,\ \beta_N:=s_{i_1}\cdots s_{i_{N-1}}(\alpha_{i_N})\,\}.
\end{equation*}
Since $\pi(B_{\bf i})=B(\infty)$ and $b_{\bf i}  : \Z_+^N \rightarrow B_{\bi}$ is a bijection by \eqref{eq:PBW vector}, one can define a crystal structure on $\Z^N_+$ by
\begin{equation}\label{eq:crystal B(infty)}
\begin{split}
& \text{ $\tf_i{\bf c}  ={\bf c}'$  if and only if $\tf_i b_{\bf i}({\bf c})\equiv b_{\bf i}({\bf c}') \!\!\mod qL(\infty)$ for ${\bf c}, {\bf c}'\in  \Z_+^N$ and $i\in I$},
\end{split}
\end{equation}
with ${\rm wt}({\bf c}) =-(c_1\beta_1 +\cdots + c_N\beta_N)$, for ${\bf c}=(c_k)\in \Z_+^N$.
We call the crystal $\Z_+^N$ the {\it crystal of ${\bf i}$-Lusztig data}, and denote it by $\cB_{\bf i}$.
{Recall that \cite{Lu93} 
\begin{equation}\label{eq:f_i for the 1st component}
\tf_{i_1}{\bf c} = (c_1+1,c_2,\cdots,c_N),\quad \text{for ${\bf c}=(c_k)\in \cB_{\bi}$}.
\end{equation}}

Let $\Omega$ be a Dynkin quiver of type $A_{n-1}$. We call a vertex $i\in I$ a sink  (resp. source) of $\Omega$ if there is no arrow going out of $i$ (resp. coming into $i$). For $i\in I$, let $s_i\Omega$ be the quiver given by reversing the arrows which end or start at $i$. 
We say that  ${\bf i}=(i_1,\ldots,i_N)\in R(w_0)$ is adapted to $\Omega$ if $i_1$ is a sink of $\Omega$, and $i_k$ is a sink of $s_{i_{k-1}}\cdots s_{i_2}s_{i_1}\Omega$ for $2\leq k\leq N$. 

Let $\cB_{\Omega}$ be the crystal $\cB_{\bf i}$ for ${\bf i}\in R(w_0)$ which is adapted to $\Omega$. Note that $\cB_{\Omega}$ is independent of the choice of ${\bf i}$ \cite{Lu90}. For ${\bf c}=(c_k)\in \cB_{\bf i}$, we write $c_{ij}=c_k$  if $\beta_k=\epsilon_i-\epsilon_j$ for $1\leq i<j\leq N$. 
{For ${\bf c}=(c_{ij})$ and ${\bf c}'=(c'_{ij})\in \cB_\Omega$, put ${\bf c}\pm{\bf c}'=(c_{ij}\pm c_{ij}')$. 
For $1\leq k<l\leq n$, let ${\bf 1}_{kl}=(c^{kl}_{ij})\in\cB_\Omega$ be such that $c^{kl}_{ij}=\delta_{ik}\delta_{jl}$.}

In the next subsections, we consider some special cases of $\Omega$, which give simple descriptions of the crystal $\cB_{\Omega}$.
\vskip 2mm

\subsection{}
We first consider the quiver $\Omega$ where all the arrows are of the same direction.
Suppose that $\Omega=\Omega^+$, where
\begin{equation*}\label{eq:orientation<-}
\Omega^+\quad : \quad
\xymatrixcolsep{2pc}\xymatrixrowsep{0pc}\xymatrix{
\bullet  & \ar@{->}[l] \bullet & \ar@{->}[l]    \cdots    & \ar@{->}[l] \bullet \\
{}_1  & {}_2 & & {}_{n-1}}.
\end{equation*}
For example, ${\bf i}=(1,2,1,3,2,1,\ldots,n-1,n-2,\ldots,2,1)$
is adapted to $\Omega^+$. We assume that $\A=\bB=[n]$ and define an injective map
\begin{equation}\label{eq:M^+(c)}
\xymatrixcolsep{3pc}\xymatrixrowsep{0pc}\xymatrix{
\cB_{\Omega^+}  \ \ar@{^{(}->}[r] & \ \M_{\A\times\bB}\ , \\
{\bf c} \ \ar@{|->}[r] & M^+({\bf c})}
\end{equation}
where $M^+({\bf c})=(m^+_{ij})$ is a strictly upper triangular matrix given by 
$m^+_{ij}=c_{ij}$ when $1\leq i<j\leq n$ and $0$ otherwise,
for ${\bf c}=(c_{ij})\in \cB_{\Omega^+}$. For $M\in \M_{\A\times\bB}$, let $M^+=(m^+_{ij})$ be the projection of $M=(m_{ij})$ onto the image of $\cB_{\Omega^+}$ under \eqref{eq:M^+(c)}, that is, $m^+_{ij}=m_{ij}$ for $1\leq i<j\leq n$, and $0$ otherwise.

Let us define $\te_i$ and $\tf_i$  for $i\in I$ on the image of $\cB_{\Omega^+}$ in $\M_{\A\times\bB}$ under \eqref{eq:M^+(c)}. 
Given ${\bf c}\in \cB_{\Omega^+}$, suppose that 
$M^+({\bf c})=M(\ba,\bb)$ for some $(\ba,\bb)\in \mathcal{I}_{\A\times\bB}$ under \eqref{eq:bijection I to M}.
Recall that ${\bb}$ is an element in a crystal $\W_{[n]}$.
For $i\in I$, we define  
\begin{equation}\label{eq:crystal<-}
\begin{split}
\te_iM^+({\bf c})&=
\begin{cases}
M(\ba,\te_i\bb)^+, & \text{if $\te_i\bb\neq {\bf 0}$},\\
{\bf 0}, & \text{if $\te_i\bb = {\bf 0}$},
\end{cases}\\
\tf_iM^+({\bf c})&=
\begin{cases}
M\left(\ba,\tf_i\bb\right), & \text{if $\tf_i\bb\neq {\bf 0}$},\\
M(\ba,\bb) + E_{i\,i+1}, & \text{if $\tf_i\bb = {\bf 0}$},
\end{cases}
\end{split}
\end{equation}
where $E_{i\,i+1}$ is an elementary matrix in $\M_{\A\times\bB}$. 
Note that ${\bf 0}$ is a formal symbol, not the zero matrix $O$.
\vskip 2mm

Next suppose that $\Omega=\Omega^-$, where 
\begin{equation*}\label{eq:orientation->}
\Omega^-\quad : \quad
\xymatrixcolsep{2pc}\xymatrixrowsep{0pc}\xymatrix{
\bullet   \ar@{->}[r] & \bullet  \ar@{->}[r]  &  \cdots   \ar@{->}[r] & \bullet \\
{}_1  & {}_2 & & {}_{n-1}}.
\end{equation*}
In this case, we assume that $\A=[n]$ and $\bB=[\ov{n}]$, and define an injective map
\begin{equation}\label{eq:M^-(c)}
\xymatrixcolsep{3pc}\xymatrixrowsep{0pc}\xymatrix{
\cB_{\Omega^-}  \ \ar@{^{(}->}[r] & \ \M_{\A\times\bB}\ , \\
{\bf c} \ \ar@{|->}[r] & M^-({\bf c})}
\end{equation}
where $M^-({\bf c})=(m^-_{ab})$ is a strictly upper triangular matrix given by 
$m^-_{n-j+1\,\ov{i}}=c_{ij}$ when $1\leq i<j\leq n$, and $0$ otherwise,
for ${\bf c}=(c_{ij})\in \cB_{\Omega^-}$.
For $M=(m_{i\ov{j}})\in \M_{\A\times\bB}$, let $M^-=(m^-_{i\ov{j}})$ be the projection of $M$ onto the image of $\cB_{\Omega^-}$ under \eqref{eq:M^-(c)}, that is, $m^-_{i\ov{j}}=m_{i\ov{j}}$ for $i+j\leq n$, and $0$ otherwise.

Let us define $\te_i$ and $\tf_i$ for $i\in I$ on the image of $\cB_{\Omega^-}$ in $\M_{\A\times\bB}$ under \eqref{eq:M^-(c)}.
Given ${\bf c}=(c_{ij})\in \cB_{\Omega^-}$, suppose that $M^-({\bf c})=M(\ba,\bb)$ for some $(\ba,\bb)\in \mathcal{I}_{\A\times\bB}$. 
For $i\in I$, we define  
\begin{equation}\label{eq:crystal->}
\begin{split}
\te_iM^-({\bf c})&= 
\begin{cases}
M(\ba,\te_i\bb)^-, & \text{if $\te_i\bb\neq {\bf 0}$},\\
{\bf 0}, & \text{if $\te_i\bb = {\bf 0}$},
\end{cases}\\
\tf_iM^-({\bf c})&=
\begin{cases}
M\left(\ba,\tf_i\bb\right), & \text{if $\tf_i\bb\neq {\bf 0}$},\\
M(\ba,\bb) + E_{n-i\,\ov{i}}, & \text{if $\tf_i\bb = {\bf 0}$},
\end{cases}
\end{split}
\end{equation}
where $E_{n-i\,\ov{i}}$ is an elementary matrix in $\M_{\A\times\bB}$.

\vskip 2mm

\begin{prop}\label{prop:crystal structure on Omega_pm}
Suppose that $\Omega$ is either $\Omega^+$ or $\Omega^-$.
The operators $\te_i$ and $\tf_i$ for $i\in I$ on $\cB_{\Omega}$ induced from \eqref{eq:crystal<-} and \eqref{eq:crystal->} coincide with those in \eqref{eq:crystal B(infty)}, that is,
 \begin{equation*}
\te_iM^\pm({\bf c}) = M^\pm\left(\te_i{\bf c}\right) \,\quad\quad \tf_iM^\pm\left({\bf c}\right) = M^\pm\left(\tf_i{\bf c}\right),
\end{equation*}
for $i\in I$ and ${\bf c}\in \cB_\Omega$. {Here we assume that $M^\pm({\bf 0})={\bf 0}$.}
\end{prop}
\pf We will consider only the case when $\Omega=\Omega^+$, since the proof for the case when $\Omega=\Omega^-$ is similar.
Let us first recall the description of \eqref{eq:crystal B(infty)} on $\cB_{\Omega}$ in \cite[Theorem 7.1]{Re} (see also \cite[Section 4.1]{S12}).
Let ${\bf c}=(c_{ij})\in \cB_{\Omega}$ be given.  For $i\in I$, put 
{\allowdisplaybreaks
\begin{equation*}
\begin{split}
c_k^{(i)} &= \sum_{s=1}^k (c_{s\,i+1}-c_{s-1\,i}) \quad \quad (1\leq k\leq i),\\
c^{(i)} & = \max\{\,c_k^{(i)}\,|\,1\leq k\leq i\,\},\\
k_0&=\min\{\,1\leq k\leq i\,|\,c^{(i)}_k=c^{(i)}\,\},\\
k_1&=\max\{\,1\leq k\leq i\,|\,c^{(i)}_k=c^{(i)}\,\},
\end{split}
\end{equation*}}
where we assume that $c_{0 i}=0$.
Then one can compute from \cite[Theorem 7.1]{Re}
{\begin{equation}\label{eq:tf_iM}
\begin{split}
\te_i{\bf c}&=
\begin{cases}
{\bf c}  + {\bf 1}_{k_0\,i} - {\bf 1}_{k_0\,i+1}, & \text{if $c^{(i)}>0$ and $k_0<i$},\\
{\bf c} - {\bf 1}_{i\,i+1}, & \text{if $c^{(i)}>0$ and $k_0=i$},\\
{\bf 0}, & \text{if $c^{(i)}=0$},
\end{cases}
\\
\tf_i{\bf c}&=
\begin{cases}
{\bf c}  - {\bf 1}_{k_1\,i} + {\bf 1}_{k_1\,i+1}, & \text{if $k_1<i$},\\
{\bf c}  + {\bf 1}_{i\,i+1}, & \text{if $k_1=i$}.
\end{cases}
\end{split}
\end{equation}}

On the other hand, suppose that $M^+({\bf c})=M(\ba,\bb)$ for some $(\ba,\bb)\in \mathcal{I}_{\A\times\bB}$ under \eqref{eq:bijection I to M}, where ${\bb}=b_1\ldots b_r$. By definition of $(\ba,\bb)$, {the subword of $\bb$ consisting of $i$ and $i+1$ is  
\begin{equation*}
\underbrace{i+1\cdots i+1}_{c_{1\,i+1}}\ \underbrace{i\cdots i}_{c_{1\,i}}\
\underbrace{i+1\cdots i+1}_{c_{2\,i+1}}\ \underbrace{i\cdots i}_{c_{2\,i}}\ \cdots\
\underbrace{i+1\cdots i+1}_{c_{i\,i+1}}.
\end{equation*}}
By the tensor product rule of crystals \eqref{eq:tensor product of crystals} (cf. \cite[Proposition 2.1.1]{KashNaka}), it is straightforward to see that 
$c^{(i)}>0$ if and only if $\varepsilon_i({\bb})>0$ and 
$\te_i{\bb}=b_1\cdots(\te_i b_s)\cdots b_r$ for some $1\leq s\leq r$ with $b_s=i+1$, {where $b_s$ is the leftmost $i+1$ in $\underbrace{i+1\cdots i+1}_{c_{k_0\,i+1}}$.} 
This implies that $\te_iM^+({\bf c})=M^+(\te_i{\bf c})$. 
Similarly, we see that
\begin{itemize}
\item[(1)] $k_1<i$ if and only if $\varphi_i({\bb})>0$ and 
$\tf_i{\bb}=b_1\cdots(\tf_i b_s)\cdots b_r$ with $b_s=i$ for some $1\leq s\leq r$, where $b_s$ is the rightmost $i$ in $\underbrace{i \cdots i }_{c_{k_1\,i}}$,

\item[(2)] $k_1=i$ if and only if $\varphi_i({\bb})=0$ or $\tf_i {\bb}={\bf 0}$,
\end{itemize}
which implies that $\tf_iM^+({\bf c})=M^+(\tf_i{\bf c})$.
\qed
\vskip 2mm

\subsection{}\label{sec:quiver with a single sink}
Now we suppose that $\Omega$ is a quiver with a single sink, that is, 
\begin{equation*}\label{eq:orientation with single sink}
\xymatrixcolsep{2pc}\xymatrixrowsep{0pc}\xymatrix{
\bullet \ar@{->}[r] &  \cdots \ar@{->}[r]  & \bullet  & \ar@{->}[l]  \cdots & \ar@{->}[l]  \bullet \\
{}_1 &  & {}_r  &  & {}_{n-1}}\ ,
\end{equation*}
for some $r\in I$. Note that we have $\Omega = \Omega^+$ if $r=1$ and $\Omega=\Omega^-$ when $r=n-1$. So we assume that $r\in I\setminus\{1,n-1\}$.
Put
\begin{equation*}\label{eq:J}
\begin{split}
&J=I\setminus\{r\},\quad
J_1=\{\,j\in J\,|\,j<r\,\},\quad J_2=\{\,j\in J\,|\,j>r\,\},
\end{split}
\end{equation*}
where $J=J_1\sqcup J_2$.
Then we have $\Phi^+=\Phi^+(J)\sqcup \Phi^+_{J_1}\sqcup \Phi^+_{J_2}$ where
\begin{equation*}\label{eq:Phi}
\begin{split}
&\Phi^+_{J_1}=\{\,\epsilon_i-\epsilon_j \,|\,1\leq  i<j \leq r\,\},\\
&\Phi^+_{J_2}=\{\,\epsilon_i-\epsilon_j \,|\,r < i<j\leq n\,\},\\
&\Phi^{+}(J)=\{\,\epsilon_i-\epsilon_j \,|\,1\leq i\leq r<j \leq n\,\}.
\end{split}
\end{equation*}
We set
\begin{equation*}
\begin{split}
\cB^J_{\Omega}&=\left\{\,{\bf c}=(c_{ij})\in \cB_\Omega\,\big\vert \,c_{ij}=0 \text{ unless $\epsilon_i-\epsilon_j\in \Phi^{+}(J)$}\,\right\},
\end{split}
\end{equation*}
which is a subcrystal of $\cB_\Omega$. {Note that we have $\tilde{f}_i{\bf c}={\bf 0}$ on $\cB^J_\Omega$ if $\tilde{f}_i{\bf c}\not\in \cB^J_\Omega$ for $i\in I$ and ${\bf c}\in \cB^J_\Omega$.}
Let $\Omega_{J_k}$ be the quiver corresponding to the vertices $J_k$ ($k=1,2$) in $\Omega$.
Then $\cB_{\Omega_{J_k}}$ is the crystal of the negative part of the quantum group $U_q(\g_{J_k})$, 
whose crystal structure is described in \eqref{eq:crystal<-} and \eqref{eq:crystal->}, respectively. 
We identify $\cB_{\Omega_{J_k}}$ with the subset of $\cB_\Omega$ consisting of ${\bf c}=(c_{ij})$ with $c_{ij}=0$ for $\epsilon_i-\epsilon_j\not\in\Phi^+_{J_k}$, and then regard it as a subcrystal of $\cB_\Omega$ where $\te_i{\bf c}=\tf_i{\bf c}={\bf 0}$ with  
$\varepsilon_i({\bf c})=\varphi_i({\bf c})=-\infty$ for $i\in J\setminus J_k$.

We define a bijection
\begin{equation}\label{eq:M(c)}
\xymatrixcolsep{3pc}\xymatrixrowsep{0pc}\xymatrix{
\cB^{J}_\Omega  \ \ar@{->}[r] & \ \M_{[\ov{r}]\times [n]\setminus [r]}\ , \\
{\bf c}  \ar@{|->}[r] & M({\bf c}) }
\end{equation}
where $M({\bf c})=(m_{ab})$ is given by 
$m_{\ov{i}\,j} = c_{ij}$ 
for ${\bf c}=(c_{ij})\in \cB_\Omega^J$ $(1\leq i\leq r<j \leq n)$. 

Given ${\bf c}\in \cB_\Omega^J$, 
suppose that 
$M({\bf c})=M(\ba,\bb)=M\left[\ba^\tau,\bb^\tau\right]=(m_{ab})$
for some $(\ba,\bb)\in \mathcal{I}_{\A\times\bB}$ (see \eqref{eq:Mtau}). 
For $i\in I$, we define
\begin{equation}\label{eq:crystal><}
\begin{split}
\te_iM({\bf c})&= 
\begin{cases}
M[\te_i\ba^\tau,\bb^\tau], & \text{if $i\in J_1$ and $\te_i\ba^\tau\neq {\bf 0}$},\\
M(\ba,\te_i\bb), & \text{if $i\in J_2$ and $\te_i\bb\neq {\bf 0}$},\\
M(\ba,\bb) - E_{\ov{r}\,r+1}, & \text{if $i=r$ and $m_{\ov{r}\,r+1}>0$},\\
{\bf 0}, & \text{otherwise},
\end{cases}\\
\tf_iM({\bf c})&=
\begin{cases}
M\left[\tf_i\ba^\tau,\bb^\tau\right], & \text{if $i\in J_1$ and $\tf_i\ba^\tau\neq {\bf 0}$},\\
M\left(\ba,\tf_i\bb\right), & \text{if $i\in J_2$ and $\tf_i\bb\neq {\bf 0}$},\\
M(\ba,\bb) + E_{\ov{r}\,r+1}, & \text{if $i=r$},\\
{\bf 0}, & \text{otherwise}.
\end{cases}
\end{split}
\end{equation}

For ${\bf c}\in\cB_\Omega$, let ${\bf c}^J$ and ${\bf c}_{J_k}$ be the restrictions of ${\bf c}$ to $\cB^J_\Omega$ and $\cB_{\Omega_{J_k}}$ $(k=1,2)$, respectively.
Then we have the following decomposition of $\cB_\Omega$ as a tensor product of its subcrystals.

\begin{thm}\label{thm:parabolic decomp}
The map 
\begin{equation*}
\xymatrixcolsep{3pc}\xymatrixrowsep{0pc}\xymatrix{
\cB_\Omega  \ar@{->}[r]  & \ \cB_\Omega^{J}\,\otimes \cB_{\Omega_{J_1}}\!\otimes \cB_{\Omega_{J_2}} \\
{\bf c}  \ar@{|->}[r] &   {\bf c}^J\otimes{\bf c}_{J_1}\otimes{\bf c}_{J_2}}
\end{equation*}
is an  isomorphism of crystals.
{Moreover, for $i\in I$ and ${\bf c}\in \cB^J_\Omega$ such that $\te_i{\bf c}\in \cB^J_\Omega$ and $\tf_i{\bf c}\in \cB^J_\Omega$, the operators $\te_i$ and $\tf_i$ on $\cB^J_{\Omega}$ induced from \eqref{eq:crystal><} coincide with those in \eqref{eq:crystal B(infty)} respectively, that is,
\begin{equation*}
\te_iM\left({\bf c}\right) = M\left(\te_i{\bf c}\right), \quad 
\tf_iM\left({\bf c}\right) = M\left(\tf_i{\bf c}\right).
\end{equation*}} 
\end{thm}
\pf 
As in Proposition \ref{prop:crystal structure on Omega_pm}, it is done by comparing with the description of \eqref{eq:crystal B(infty)} on $\cB_{\Omega}$ using \cite[Theorem 7.1]{Re}. 

For ${\bf c}\in \cB_\Omega$,
let $\psi({\bf c})={\bf c}^J\otimes{\bf c}_{J_1}\otimes{\bf c}_{J_2}$. It is clear that $\psi$ is a bijection. So it remains to show that $\psi$ commutes with $\te_i$ and $\tf_i$ for $i\in I$.

Suppose that ${\bf c}=(c_{ij})\in \cB_\Omega$ is given. 
First, {we have by \eqref{eq:f_i for the 1st component}} 
\begin{equation*}\label{eq:}
\begin{split}
\te_r{\bf c}&=
\begin{cases}
{\bf c}  - {\bf 1}_{r\,r+1}, & \text{if $c_{r\,r+1}>0$},\\
{\bf 0}, & \text{if $c_{r\,r+1}=0$},
\end{cases}\quad \quad
\tf_r{\bf c}= 
{\bf c}  + {\bf 1}_{r\,r+1},
\end{split}
\end{equation*}
which immediately implies that
\begin{equation*}\label{eq:M commutes with te_r and tf_r}
\te_r M\left( {\bf c}^J \right) = M\left(\te_r{\bf c}^J\right),\quad 
\tf_rM\left({\bf c}^J\right) = M\left(\tf_r{\bf c}^J\right), 
\end{equation*}
assuming that $M({\bf 0})={\bf 0}$, and hence $\psi$ commutes with $\te_r$ and $\tf_r$.

Next, we fix $i\in J_1$. Let
\begin{equation*}\label{eq:crytal on B_Omega-1}
\begin{split}
c_k^{(i)}&= 
\begin{cases}
c_{i\,r+1}, & \text{if $k=1$},\\
c_1^{(i)}+\sum_{s=2}^{k} (c_{i\, r+s}-c_{i+1\, r+s-1}), & \text{if $2\leq k\leq n-r$},\\
c^{(i)}_{n-r}+ (c_{i\,r}-c_{i+1\,n}), & \text{if $k=n-r+1$},\\
c^{(i)}_{n-r+1}+\sum_{s=1}^{k-n+r-1} (c_{i\,r-s}-c_{i+1\,r-s+1}),\!\! & \text{if $n-r+2\leq k\leq n-i$},
\end{cases}
\end{split}
\end{equation*}
and
{\allowdisplaybreaks
\begin{equation}\label{eq:crytal on B_Omega-2}
\begin{split}
c^{(i)}& = \max\{\,c_k^{(i)}\,|\,1\leq k\leq n-i\,\},\\
k_0&=\min\{\,1\leq k\leq n-i\,|\,c^{(i)}_k=c^{(i)}\,\},\\
k_1&=\max\{\,1\leq k\leq n-i\,|\,c^{(i)}_k=c^{(i)}\,\}.
\end{split}
\end{equation}
Note that if $c^{(i)}>0$, then we have $c_{i\,k_0+r}>0$ when $k_0\leq n-r$, and $c_{i\,n-k_0+1}>0$ when $k>n-r$. 
Also if $k_0>n-r$, then we necessarily have $c^{(i)}>0$.
By \cite[Theorem 7.1]{Re}, one can compute directly that
\begin{equation}\label{eq:crytal on B_Omega-3}
\begin{split}
\te_i{\bf c}&=
\begin{cases}
{\bf c}  - {\bf 1}_{i\,k_0+r} + {\bf 1}_{i+1\,k_0+r}, & \text{if $c^{(i)}>0$ and $k_0\leq n-r$},\\
{\bf c}  - {\bf 1}_{i\,n-k_0+1} + {\bf 1}_{i+1\,n-k_0+1}, & \text{if  $n-r+1\leq k_0 \leq n-i-1$},\\
{\bf c}  - {\bf 1}_{i\,i+1}, & \text{if $k_0=n-i$},\\
{\bf 0}, & \text{if $c^{(i)}=0$},
\end{cases}\\
\tf_i{\bf c}&=
\begin{cases}
{\bf c}  + {\bf 1}_{i\,k_1+r} - {\bf 1}_{i+1\,k_1+r}, & \text{if $k_1\leq n-r$},\\
{\bf c}  + {\bf 1}_{i\,n-k_1+1} - {\bf 1}_{i+1\,n-k_1+1}, & \text{if $n-r+1\leq k_1 \leq n-i-1$},\\
{\bf c}  + {\bf 1}_{i\,i+1}, & \text{if $k_1=n-i$,}
\end{cases}\\\end{split}
\end{equation}
\noindent (see for example, the Auslander-Reiten quiver in Example \ref{ex:example of embdedding}, which might be helpful to see which $c_{ij}$'s are involved for $\te_i$ and $\tf_i$, and how they are arranged).\vskip 2mm}

{\it Case 1}. Suppose that ${\bf c}={\bf c}^J\in \cB^J_\Omega$, that is, $c_{ij}=0$ unless $\epsilon_i-\epsilon_j \in \Phi^+(J)$. We have 
\begin{equation*}
c^{(i)}_{n-r}\geq c^{(i)}_{n-r+1}=\cdots=c^{(i)}_{n-i},
\end{equation*}
which implies that $k_0\leq n-r$. Note that if $k_1>n-r$, then we have $c_{i+1\,n}=0$ and hence $k_1=n-i$. 

Let $(\ba,\bb)\in \mathcal{I}_{[\ov{r}]\times ([n]\setminus [r])}$ be such that $M({\bf c})=M[\ba^\tau,\bb^\tau]$. {Note that the subword of ${\ba}^\tau$ consisting of $i$ and $i+1$ is
\begin{equation*}
\begin{split}
&\underbrace{\ov{i}\cdots \ov{i}}_{c_{i,r+1}}\ \underbrace{\ov{i+1}\cdots \ov{i+1}}_{c_{i+1,r+1}}\ \cdots\
\underbrace{\ov{i}\cdots \ov{i}}_{c_{i\,n}}\ \underbrace{\ov{i+1}\cdots \ov{i+1}}_{c_{i+1\,n}}.  
\end{split}
\end{equation*}}
By the tensor product rule \eqref{eq:tensor product of crystals}, 
we have $\varepsilon_i({\bf a}^\tau)=c^{(i)}$ and 
\begin{equation*}\label{eq:M commutes with te}
M(\te_i{\bf c})=M\left[\te_i\ba^\tau,{\bb}^\tau\right]=\te_iM({\bf c}),
\end{equation*}
If $k_1\leq n-r$, then  $\tf_i{\bf a}^\tau\neq {\bf 0}$ and
\begin{equation*}\label{eq:M commutes with tf}
M\left(\tf_i{\bf c}\right)=M\left[\tf_i\ba^\tau,{\bb}^\tau\right]=\tf_iM({\bf c}).
\end{equation*}
If $k_1>n-r$, then $\tf_i{\bf a}^\tau={\bf 0}$ and $\tf_i {\bf c}\not\in \cB^J_{\Omega}$, which implies that $\tf_iM({\bf c})={\bf 0}$ and $\tf_i{\bf c}={\bf 0}$ in $\cB^J_\Omega$, respectively. We have $\tf_iM({\bf c})=M\left(\tf_i{\bf c}\right)={\bf 0}$. 

\vskip 2mm

{\it Case 2}. Suppose that ${\bf c}\in \cB_\Omega$ is arbitrary. 
We assume that 
\begin{equation*}
M_1:=M\left({\bf c}^J\right)=M[\ba^\tau,\bb^\tau],\quad M_2:=M^-({\bf c}_{J_1})=M({\bf a}',{\bf b}'),
\end{equation*}
for some $(\ba,\bb)\in \mathcal{I}_{[\ov{r}]\times ([n]\setminus [r])}$ 
and $(\ba',\bb')\in \mathcal{I}_{[r]\times [\ov{r}]}$. 
By Proposition \ref{prop:crystal structure on Omega_pm} and the arguments in {\it Case 1}, 
we see that  
\begin{equation}\label{eq:varepsilon=}
\varepsilon_i(M_1)=\varepsilon_i(\ba^\tau),\quad
\varepsilon_i(M_2)=\varepsilon_i(\bb').
\end{equation}
Since $\langle {\rm wt}(M_1),h_i \rangle =\langle {\rm wt}(\ba^\tau),h_i \rangle$ and 
$\langle {\rm wt}(M_2),h_i \rangle =\langle {\rm wt}(\bb'),h_i \rangle$, we also have
\begin{equation}\label{eq:varphi=}
\varphi_i(M_1)=\varphi_i(\ba^\tau),\quad \varphi_i(M_2)=\varphi_i(\bb').
\end{equation}
By \eqref{eq:crytal on B_Omega-2} and \eqref{eq:crytal on B_Omega-3}, we have
\begin{equation*}
\begin{split}
\psi(\te_i{\bf c}) = \left(\te_i{\bf c}^J\right)\otimes {\bf c}_{J_1} \otimes {\bf c}_{J_2}
\quad & \Longleftrightarrow \quad k_0\leq n-r   \\
& \Longleftrightarrow \quad \te_i(\ba^\tau\otimes \bb') =(\te_i\ba^\tau)\otimes \bb' \\
& \Longleftrightarrow \quad \varphi_i(\ba^\tau)\geq \varepsilon_i(\bb'), \\
\end{split}
\end{equation*}
\begin{equation*}
\begin{split}
\psi(\te_i{\bf c}) = {\bf c}^J \otimes (\te_i{\bf c}_{J_1}) \otimes {\bf c}_{J_2}
\quad & \Longleftrightarrow \quad k_0> n-r   \\
& \Longleftrightarrow \quad \te_i(\ba^\tau\otimes \bb') =\ba^\tau \otimes(\te_i \bb') \\
& \Longleftrightarrow \quad \varphi_i(\ba^\tau) < \varepsilon_i(\bb'). \\
\end{split}
\end{equation*}
Therefore, we have by \eqref{eq:varepsilon=} and \eqref{eq:varphi=}
\begin{equation*}
\psi\left(\te_i{\bf c}\right)=
\begin{cases}
(\te_i{\bf c}^J)\otimes {\bf c}_{J_1} \otimes {\bf c}_{J_2},& 
\text{if $\varphi_i(M_1) \geq \varepsilon_i(M_2)$}, \\
{\bf c}^J \otimes (\te_i{\bf c}_{J_1}) \otimes {\bf c}_{J_2},& 
\text{if $\varphi_i(M_1) < \varepsilon_i(M_2)$}.
\end{cases}
\end{equation*}
Similarly, we have
\begin{equation*}
\psi\left(\tf_i{\bf c}\right)=
\begin{cases}
\left(\tf_i{\bf c}^J\right)\otimes {\bf c}_{J_1} \otimes {\bf c}_{J_2},& 
\text{if $\varphi_i(M_1) > \varepsilon_i(M_2)$}, \\
{\bf c}^J \otimes \left(\tf_i{\bf c}_{J_1}\right) \otimes {\bf c}_{J_2},& 
\text{if $\varphi_i(M_1) \leq  \varepsilon_i(M_2)$}.
\end{cases}
\end{equation*}
It follows that $\psi$ commutes with $\te_i$ and $\tf_i$ for $i\in J_1$.

By the same arguments, we can show that $\psi$ commutes with $\te_i$ and $\tf_i$ for $i\in J_2$. This completes the proof. \qed\vskip 2mm

 


Let $B_J(\infty)$ denote the $\g_J$-crystal of the negative part of $U_q(\g_J)$, and {extend it to a $\g$-crystal with $\te_i b= \tf_i b={\bf 0}$ and $\varepsilon_i(b)=\varphi_i(b)=-\infty$ for $i\in I\setminus J$ and $b\in B_J(\infty)$.} 

Let $W_J$ be the Weyl group of $\g_J$ generated by $s_j$ for $j\in J$, and let $w^J$ be the longest element in the set of coset representatives of minimal length in $W/W_J$. 
Consider $\bi\in R(w_0)$ corresponding to $w_0=w^Jw_{J}$, where $w_{J}$ is the longest element in $W_J$. Let $U^-\left(J\right)$ be the $\mathbb{Q}(q)$-subspace of $U^-$ spanned by $b_{\bi}({\bf c})\in B_{\bi}$ for ${\bf c}\in \Z_+^N$ such that ${\bf c}={\bf c}^J$. Then $U^-\left(J\right)$ is independent of the choice of $\bi$, and forms a subalgebra of $U^-$ called the {\it quantum nilpotent subalgebra} associated to $w^J$ \cite{Lu93}.
By using a PBW basis, we see that the multiplication in $U^-$ gives an isomorphism 
of a $\mathbb{Q}(q)$-vector space
\begin{equation}\label{eq:U^-_decomp}
U^- \cong U^-\left(J\right)\otimes U^-_J
\end{equation}
{(see \cite{Ki,T} for more details and its generalization to the case of a symmetrizable Kac-Moody algebra).} 
The image of a PBW basis of $U^-(J)$ under the canonical projection $\pi$ 
forms a subcrystal of $\pi(B_{\bi})=B(\infty)$ in $L(\infty)/qL(\infty)$, which we denote by $B^J(\infty)$. Then we have the following tensor product decomposition of $B(\infty)$, which is a crystal version of \eqref{eq:U^-_decomp}.

\begin{cor}
As a $\g$-crystal, we have
\begin{equation*}
B(\infty)\cong B^J(\infty)\otimes B_J(\infty).
\end{equation*}
\end{cor}
\pf We have 
$B_J(\infty)\cong B_{J_1}(\infty)\otimes B_{J_2}(\infty)\cong 
\cB_{\Omega_{J_1}}\!\otimes \cB_{\Omega_{J_2}}$ and $B^J(\infty)\cong \cB_{\Omega}^J$.
Hence it follows from Theorem \ref{thm:parabolic decomp}.
\qed 

\begin{rem}{\rm
The isomorphism in Theorem \ref{thm:parabolic decomp} is a special case of the bijection $\Omega_w$ for $w\in W$ in \cite[Proposition 5.25]{BKT} when $w=w^J$ (see also \cite[Proposition 3.14]{Ki}). We should remark that $\Omega_w$ is not in general a crystal isomorphism for arbitrary $w\in W$. {For example, suppose that $n=3$ and $\bi=(1,2,1)$, that is, $\beta_1=\alpha_1$, $\beta_2=\alpha_1+\alpha_2$, $\beta_3=\alpha_2$.  If $w=s_1$, then $\Omega_w$ is given by sending ${\bf c}$ to $({\bf c}_{\leq w}, {\bf c}_{> w})$ where ${\bf c}_{\leq w}=(c_1,0,0)$ and ${\bf c}_{> w}=(0,c_2,c_3)$ for ${\bf c}=(c_1,c_2,c_3)$. But the mapping ${\bf c} \mapsto {\bf c}_{\leq w}\otimes {\bf c}_{>w}$ does not define a morphism of crystals in this case (see Proposition \ref{prop:crystal structure on Omega_pm}). It would be interesting to characterize $w\in W$ such that the map $\Omega_w : {\bf c} \mapsto {\bf c}_{\leq w}\otimes {\bf c}_{>w}$ is an isomorphism of crystals by using the result in \cite{SST}, where a connection between crystal structure of Lusztig data and signature rule in tensor product is studied.}}
\end{rem}

\section{Crystal embedding of Young tableaux into Lusztig data}\label{sec:embedding}

\subsection{}

Let $\la\in \cP_n$ be given. For $S\in SST_{[n]}(\la)$, we define  
\begin{equation}\label{eq:c^+}
{\bf c}^+(S)=(c_{ij})\in \cB_{\Omega^+},
\end{equation}
where $c_{ij}$ is given by the number of $j$'s appearing in the $i$th row of $S$ for $1\leq i<j\leq n$. 
Then we have the following, which is already well-known to experts in this area and which the author learned from Y. Saito. 

\begin{prop}\label{prop:embedding Omega^+} For $\la\in \cP_n$, the map 
\begin{equation*}
\xymatrixcolsep{3pc}\xymatrixrowsep{0pc}\xymatrix{
SST_{[n]}(\la)\otimes T_{-\la}  \ \ar@{->}[r] & \ \cB_{\Omega^{+}}\ , \\
S\otimes t_{-\la} \ar@{|->}[r] & {\bf c}^+(S)}
\end{equation*}
is an embedding of crystals.
\end{prop}
\pf It follows immediately from comparing the crystal structures on $SST_{[n]}(\la)$ and $\cB_{\Omega^+}$ described in Proposition \ref{prop:crystal structure on Omega_pm}.
\qed\vskip 2mm

We also have an embedding into $\cB_{\Omega^-}$.
For  $T\in SST_{[\ov{n}]}(\la)$, we define
\begin{equation}\label{eq:c^-}
{\bf c}^-(T)=(c_{ij})\in \cB_{\Omega^-},
\end{equation}
where $c_{ij}$ is given by the number of $\ov{i}$'s appearing in the $(n-j+1)$th row of $T$ for $1\leq i<j\leq n$. Similarly, for $S\in SST_{[n]}(\la)$, we define
\begin{equation}\label{eq:c_-}
{\bf c}_-(S)={\bf c}^-\left(\sigma^{-d}(S)^{\nw}\right)\in \cB_{\Omega^-},
\end{equation}
for some $d\geq \lambda_1$. Note that ${\bf c}_-(S)$ does not depend on the choice of $d$.

\begin{prop}\label{prop:embedding Omega^-}
For $\la\in \cP_n$, the maps
\begin{equation*}
\xymatrixcolsep{3pc}\xymatrixrowsep{0pc}\xymatrix{
SST_{[\ov{n}]}(\la)\otimes T_{w_0\la} \ \ar@{->}[r] & \ \cB_{\Omega^{-}}\ , \\
T\otimes t_{w_0\la} \ar@{|->}[r] & {\bf c}^-(T)}
\quad\quad
\xymatrixcolsep{3pc}\xymatrixrowsep{0pc}\xymatrix{
SST_{[n]}(\la)\otimes T_{-\la} \ \ar@{->}[r] & \ \cB_{\Omega^{-}} \\
S\otimes t_{-\la} \ar@{|->}[r] & {\bf c}_-(S)}
\end{equation*}
are embeddings of crystals.
\end{prop}
\pf  It follows from \eqref{eq:tau} and Proposition \ref{prop:crystal structure on Omega_pm}.\qed 
 
\begin{ex}\label{ex:embedding-1}
{\rm
Suppose that $n=6$ and let\vskip 2mm
$$ S\ =\
\resizebox{.18\hsize}{!}{$
{\def\lr#1{\multicolumn{1}{|@{\hspace{.75ex}}c@{\hspace{.75ex}}|}{\raisebox{-.04ex}{$#1$}}}\raisebox{-.6ex}
{$\begin{array}{cccccc}
\cline{1-6}
\lr{1} &\lr{1} & \lr{1} &\lr{2} &\lr{2} &\lr{3} \\
\cline{1-6}
\lr{2} & \lr{3} & \lr{3} & \lr{5} & \lr{6} \\
\cline{1-5}
\lr{4} & \lr{4} & \lr{4} \\
\cline{1-3}
\lr{5}& \lr{5} & \lr{6} \\
\cline{1-3}
\lr{6} & \lr{6}  \\
\cline{1-2} 
\end{array}$}}$}\ \in SST_{[6]}(6,5,3,3,2)\ .$$\vskip 2mm

Then we have by \eqref{eq:c^+}

\begin{equation*}
{\bf c}^+(S) = 
\begin{bmatrix}
c_{12} & c_{13} & c_{14} & c_{15} & c_{16} \\
       & c_{23} & c_{24} & c_{25} & c_{26} \\
       &        & c_{34} & c_{35} & c_{36} \\
       &        &        & c_{45} & c_{46} \\
       &        &        &        & c_{56} \\
\end{bmatrix}
=
\begin{bmatrix}
2 & 1 & 0 & 0 & 0 \\
  & 2 & 0 & 1 & 1 \\
  &   & 3 & 0 & 0 \\
  &   &   & 2 & 1 \\ 
  &   &   &   & 2 \\  
\end{bmatrix}\ \in \cB_{\Omega^+}\ .
\end{equation*}\vskip 2mm

\noindent On the other hand,
$$ \sigma^{-6}(S)\ =\
\resizebox{.17\hsize}{!}{$
{\def\lr#1{\multicolumn{1}{|@{\hspace{.75ex}}c@{\hspace{.75ex}}|}{\raisebox{-.04ex}{$#1$}}}\raisebox{-.6ex}{$\begin{array}{cccccc}
\cline{6-6}
& & & & & \lr{\ov{6}} \\
\cline{4-6}
& & & \lr{\ov{6}} & \lr{\ov{5}} & \lr{\ov{5}} \\
\cline{4-6}
& & & \lr{\ov{4}} & \lr{\ov{4}} & \lr{\ov{4}} \\
\cline{3-6}
& & \lr{\ov{5}} & \lr{\ov{3}} & \lr{\ov{3}} & \lr{\ov{2}} \\
\cline{1-6}
\lr{\ov{3}} & \lr{\ov{2}} & \lr{\ov{2}} & \lr{\ov{1}} & \lr{\ov{1}} & \lr{\ov{1}} \\
\cline{1-6}
\end{array}$}}$}\ \ , \ \ \ \ \ \
\sigma^{-6}(S)^{\nw}\ = \
\resizebox{.17\hsize}{!}{$
{\def\lr#1{\multicolumn{1}{|@{\hspace{.75ex}}c@{\hspace{.75ex}}|}{\raisebox{-.04ex}{$#1$}}}\raisebox{-.6ex}{$\begin{array}{cccccc}
\cline{1-6}
\lr{\ov{6}} & \lr{\ov{6}} & \lr{\ov{5}} & \lr{\ov{4}} & \lr{\ov{1}} & \lr{\ov{1}}  \\
\cline{1-6}
\lr{\ov{5}} & \lr{\ov{5}} & \lr{\ov{4}} & \lr{\ov{2}}  \\
\cline{1-4}
\lr{\ov{4}} & \lr{\ov{3}} & \lr{\ov{3}} \\
\cline{1-3}
\lr{\ov{3}} & \lr{\ov{2}} & \lr{\ov{1}} \\
\cline{1-3}
\lr{\ov{2}}  \\
\cline{1-1}
\end{array}$}}$}\ .
$$\vskip 2mm
\noindent  Hence by \eqref{eq:c_-}, we have
\begin{equation*}
{\bf c}_-(S) = 
\begin{bmatrix}
c_{56} & c_{46} & c_{36} & c_{26} & c_{16} \\
       & c_{45} & c_{35} & c_{25} & c_{15} \\
       &        & c_{34} & c_{24} & c_{14} \\
       &        &        & c_{23} & c_{13} \\
       &        &        &        & c_{12} \\
\end{bmatrix}
=
\begin{bmatrix}
1 & 1 & 0 & 0 & 2 \\
  & 1 & 0 & 1 & 0 \\
  &   & 2 & 0 & 0 \\
  &   &   & 1 & 1 \\ 
  &   &   &   & 0 \\  
\end{bmatrix}\ \in \cB_{\Omega^-}\ .
\end{equation*}\vskip 2mm

}
\end{ex}

\subsection{}
Let $\Omega$ be a quiver with a single sink 
\begin{equation*}
\xymatrixcolsep{2pc}\xymatrixrowsep{0pc}\xymatrix{
\bullet \ar@{->}[r] &  \cdots \ar@{->}[r]  & \bullet  & \ar@{->}[l]  \cdots & \ar@{->}[l]  \bullet \\
{}_1 &  & {}_r  &  & {}_{n-1}}
\end{equation*}
for some $r\in I\setminus\{1,n-1\}$. We keep the notations in Section \ref{sec:quiver with a single sink}. 

Let $\la=(\la_1,\ldots,\la_n)\in \cP_n$ be given. Choose $d\geq \la_1$ and put
\begin{equation}\label{eq:eta zeta for la}
\eta=(d-\la_{r},\ldots,d-\la_1)\in \cP_r,\quad \zeta=(\la_{r+1},\ldots,\la_n)\in \cP_{n-r}.
\end{equation}
We define a map
\begin{equation}\label{eq:embedding_1}
\xymatrixcolsep{3pc}\xymatrixrowsep{0pc}\xymatrix{
SST_{[n]}(\la)  \ \ar@{->}[r] & \ SST_{[n]\setminus[r]}(\zeta)\times SST_{[\ov{r}]}(\eta)\times \M_{[\ov{r}]\times ([n]\setminus[r])} \\
S \ \ar@{|->}[r] & \  (S^+,S^-,M) }
\end{equation}
where $(S^+,S^-,M)$ is determined by the following steps:
\begin{itemize}
\item[(i)] let $S^+\in \ SST_{[n]\setminus[r]}(\zeta)$ be given by removing the first $r$ rows in $S$,

\item[(ii)] let $S\setminus S^+$ denote the subtableau of $S$ obtained by removing $S^+$, and put 
\begin{equation*}
\begin{split}
P' &= \text{the subtableau of $S\setminus S^+$ with entries in $[r]$},\\
Q &= \text{the subtableau of $S\setminus S^+$ with entries in $[n]\setminus[r]$},\\
\end{split}
\end{equation*}

\item[(iii)] putting $P=\sigma^{-d}(P')$ (see \eqref{eq:tau}), 
we have for some $\nu\in\cP_r$ with $\eta\subset\nu$
\begin{equation*}
\begin{split}
(P,Q)\in SST_{[\ov{r}]}(\nu^\pi) \times SST_{[n]\setminus[r]}(\left(\nu/\eta\right)^\pi),
\end{split}
\end{equation*}

\item[(iv)] applying $\kappa^{-1}$ in \eqref{eq:skew_RSK_antinormal}, we get 
\begin{equation*}
\begin{split}
(T,M)=\kappa^{-1}(P,Q) \in SST_{[\ov{r}]}(\eta^\pi)\times \M_{[\ov{r}]\times ([n]\setminus[r])},
\end{split}
\end{equation*}

\item[(v)] let $S^-= T^{\nw} \in SST_{[\ov{r}]}(\eta)$.
\end{itemize}\vskip 2mm

It can be summarized as follows:
\begin{equation}\label{eq:embedding}
\xymatrixcolsep{3pc}\xymatrixrowsep{3pc}\xymatrix{
S  \  \ar@{->}[r]^{{\rm (i)}} \ar@{.>}[d]  & (S^+,S\setminus S^+) \  \ar@{->}[r]^{{\rm (ii)}} &  (S^+,P',Q)  
\ar@{->}[d]^{{\rm (iii)}} \\
 (S^+,S^-,M) & \ar@{->}[l]^{{\rm (v)}} (S^+,T,M)  & \ar@{->}[l]^{{\rm (iv)}} (S^+,P,Q)    }
\end{equation}
{Note that the step (i) is injective, and the other steps (ii), (iii), (iv), (v) are bijective by definition.} 
Hence  the map \eqref{eq:embedding_1} is injective, where
\begin{equation*}
{\rm wt}(S^+)+{\rm wt}(S^-)+{\rm wt}(M) = {\rm wt}(S) -d(\epsilon_1+\cdots+\epsilon_r).
\end{equation*}

Now for $S\in SST_{[n]}(\la)$ which is mapped to $(S^+,S^-,M)$ under \eqref{eq:embedding_1}, we define  
\begin{equation}\label{eq:c}
{\bf c}(S)\in \cB_{\Omega},
\end{equation}
to be the unique ${\bf c}\in \cB_\Omega$ such that
\begin{itemize}
\item[(1)] $M\left({\bf c}^J\right)=M$ under \eqref{eq:M(c)},

\item[(2)] ${\bf c}_{J_1}={\bf c}^-(S^-)$ and ${\bf c}_{J_2}={\bf c}^+(S^+)$ under \eqref{eq:c^-} and \eqref{eq:c^+}, respectively.
\end{itemize}
Note that $S^-$ depends on $d$, but ${\bf c}^-(S^-)$ does not.
Then we have the following, which is the main result in this paper.

\begin{thm}\label{thm:main-1} For $\la\in\cP_n$, the map
\begin{equation*}
\xymatrixcolsep{3pc}\xymatrixrowsep{0pc}\xymatrix{
SST_{[n]}(\la)\otimes T_{-\la}  \ \ar@{->}[r] & \ \cB_\Omega \\
S\otimes t_{-\la} \ar@{|->}[r] & {\bf c}(S) }
\end{equation*}
is an embedding of crystals, where ${\bf c}(S)$ is given in \eqref{eq:c}.
\end{thm}
\pf Put $\M=\M_{[\ov{r}]\times ([n]\setminus[r])}$.
Note that $\M_{[\ov{r}]\times ([n]\setminus[r])}$ has a crystal structure isomorphic to $\cB_{\Omega}^J$ induced from the bijection \eqref{eq:M(c)}, which can be described as in \eqref{eq:crystal><} by Theorem \ref{thm:parabolic decomp}.


Choose $d\geq \la_1$ and let $\eta$ and $\zeta$ be as in \eqref{eq:eta zeta for la}.
Define a $\g$-crystal
\begin{equation*}
\M_\la = \M \otimes SST_{[\ov{r}]}(\eta)\otimes SST_{[n]\setminus[r]}(\zeta),
\end{equation*}
where we extend a $\g_{J_1}$-crystal $SST_{[\ov{r}]}(\eta)$ 
and a $\g_{J_2}$-crystal $SST_{[n]\setminus[r]}(\zeta)$ to $\g$-crystals in a trivial way.
Put $\xi=-d(\epsilon_1+\cdots+\epsilon_r)$.

By \eqref{eq:crystal><}, we see that the crystal structure on $\cM_\la$ coincides with the one given in \cite[Section 4.2]{K09}. Moreover, if we put ${H}_\la={O}\otimes b_1\otimes b_2$ where ${O}$ is the zero matrix, $b_1$ (resp. $b_2$) is the highest weight element in $SST_{[\ov{r}]}(\eta)$ (resp. $SST_{[n]\setminus[r]}(\zeta)$) with weight $-\eta_r\epsilon_1-\cdots-\eta_1\epsilon_r$ (resp. $\zeta_1\epsilon_{r+1}+\cdots+\zeta_{n-r}\epsilon_{n}$), then ${H}_\la$ is the highest weight element with weight $\la+\xi$ and 
\begin{equation*}
\M_\la =\{\,\tf_{i_1}\cdots \tf_{i_r} {H}_\la\,|\,r\geq 0, i_1,\ldots,i_r\in I\,\},
\end{equation*}
\cite[Proposition 4.5]{K09}.
By \cite[Proposition 4.6]{K09}, the map
\begin{equation}\label{eq:map1}
\xymatrixcolsep{3pc}\xymatrixrowsep{0pc}\xymatrix{
SST_{[n]}(\la)\otimes T_{\xi}  \ \ar@{->}[r] & \ \M_\la \\
S\otimes t_{\xi} \ar@{|->}[r] & M\otimes S^-\otimes S^+ }
\end{equation}
is an embedding of $\g$-crystals, where $(S^+,S^-,M)$ is the triple associated to $S$ in \eqref{eq:embedding_1}. Taking tensor product by $T_{-\la-\xi}$ and then applying 
\eqref{eq:M(c)}, Propositions \ref{prop:embedding Omega^+} and \ref{prop:embedding Omega^-}, we have an embedding
\begin{equation}\label{eq:map2}
\xymatrixcolsep{3pc}\xymatrixrowsep{0pc}\xymatrix{
\M_\la \otimes T_{-\xi-\lambda}  \ \ar@{->}[r] & \ \cB_\Omega^{J}\,\otimes \cB_{\Omega_{J_1}}\!\otimes \cB_{\Omega_{J_2}}\ , \\
M\otimes S^-\otimes S^+ \otimes t_{-\xi-\lambda} \ar@{|->}[r] & {\bf c}^J\otimes{\bf c}_{J_1}\otimes{\bf c}_{J_2}}
\end{equation}
where ${\bf c}={\bf c}(S)$ is given in \eqref{eq:c}. Finally composing \eqref{eq:map1}, \eqref{eq:map2}, and then the inverse of the map in Theorem \ref{thm:parabolic decomp}, we obtain the required embedding.
\qed

\begin{rem}{\rm
Suppose that $\Omega$ is a quiver with a single source. Let $\cB_{\Omega}^\ast$ be the set $\cB_\Omega$ with the $\ast$-crystal structure \cite{Kas93}. By similar methods as in Theorem \ref{thm:main-1}, we can construct an embedding of $SST_{[n]}(\la)$ into $\cB_{\Omega}^\ast$ for $\la\in \cP_n$.   
}
\end{rem}

\begin{ex}\label{ex:example of embdedding}
{\rm

Suppose that $\Omega$ is given by
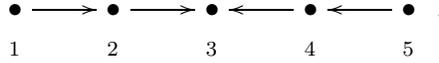
\begin{equation*}
\xymatrixcolsep{2pc}\xymatrixrowsep{0pc}\xymatrix{
\bullet \ar@{->}[r] &  \bullet \ar@{->}[r]  & \bullet  & \ar@{->}[l]  \bullet & \ar@{->}[l]  \bullet \\
{}_1 & {}_2 & {}_3  & {}_4 & {}_{5}}\ .
\end{equation*}
Recall that the Auslander-Reiten quiver of representations of $\Omega$ is  
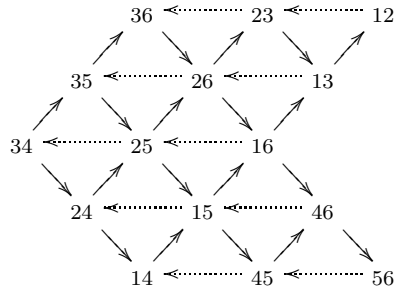
\begin{equation*} 
\xymatrixcolsep{0.5pc}\xymatrixrowsep{1pc}\xymatrix{
 &    & {}_{36} \ar@{->}[rd] &   & {}_{23} \ar@{.>}[ll] \ar@{->}[rd] & & {}_{12} \ar@{.>}[ll] \\
 &  {}_{35} \ar@{->}[ur] \ar@{->}[rd]  &    &  {}_{26} \ar@{.>}[ll] \ar@{->}[ur] \ar@{->}[rd] &  & {}_{13} \ar@{.>}[ll] \ar@{->}[ur]  \\
{}_{34} \ar@{->}[ur] \ar@{->}[rd] &    & {}_{25} \ar@{.>}[ll] \ar@{->}[ur] \ar@{->}[rd] &   &  {}_{16} \ar@{.>}[ll] \ar@{->}[ur] \ar@{->}[rd]\\
 &  {}_{24}  \ar@{->}[ur] \ar@{->}[rd] &    &   {}_{15} \ar@{.>}[ll] \ar@{->}[ur] \ar@{->}[rd] &  & {}_{46}\ar@{.>}[ll]\ar@{->}[rd]\\  
 &    & {}_{14} \ar@{->}[ur]   &   & {}_{45} \ar@{.>}[ll] \ar@{->}[ur]  & & {}_{56} \ar@{.>}[ll] }
\end{equation*} 
\vskip 2mm
\noindent which might be helpful for the reader to see \eqref{eq:crytal on B_Omega-3} from Reineke's description of $B(\infty)$ \cite{Re}. Here the vertex $``ij"$ denotes the indecomposable representation of $\Omega$ corresponding to the  positive root $\epsilon_i-\epsilon_j\in\Phi^+$ for $1\leq i<j\leq 6$, the solid arrows denote the morphisms between them, and the dotted arrows denote the Auslander-Reiten translation functor denoted by $\tau$ in \cite{Re}.

Let $S$ be as in Example \ref{ex:embedding-1}.
Let us apply the map \eqref{eq:embedding_1} to $S$ following the steps in  \eqref{eq:embedding}. First, we have
$$ S^+ \ =\
\resizebox{.09\hsize}{!}{$
{\def\lr#1{\multicolumn{1}{|@{\hspace{.75ex}}c@{\hspace{.75ex}}|}{\raisebox{-.04ex}{$#1$}}}\raisebox{-.6ex}
{$\begin{array}{ccc}
\cline{1-3}
\lr{5}& \lr{5} & \lr{6} \\
\cline{1-3}
\lr{6} & \lr{6}  \\
\cline{1-2} 
\end{array}$}}$}\ \ ,\ \ \ \ \ \  \ S\setminus S^+\ = \
\resizebox{.18\hsize}{!}{$
{\def\lr#1{\multicolumn{1}{|@{\hspace{.75ex}}c@{\hspace{.75ex}}|}{\raisebox{-.04ex}{$#1$}}}\raisebox{-.6ex}
{$\begin{array}{cccccc}
\cline{1-6}
\lr{1} &\lr{1} & \lr{1} &\lr{2} &\lr{2} &\lr{3} \\
\cline{1-6}
\lr{2} & \lr{3} & \lr{3} & \lr{5} & \lr{6} \\
\cline{1-5}
\lr{4} & \lr{4} & \lr{4} \\
\cline{1-3}
\end{array}$}}$}\ .$$\vskip 2mm
\noindent Separating $S\setminus S^+$ into subtableaux with entries in $\{1,2,3\}$ and $\{4,5,6\}$, we get
$$ P'\ = \
\resizebox{.18\hsize}{!}{$
{\def\lr#1{\multicolumn{1}{|@{\hspace{.75ex}}c@{\hspace{.75ex}}|}{\raisebox{-.04ex}{$#1$}}}\raisebox{-.6ex}
{$\begin{array}{cccccc}
\cline{1-6}
\lr{1} &\lr{1} & \lr{1} &\lr{2} &\lr{2} &\lr{3} \\
\cline{1-6}
\lr{2} & \lr{3} & \lr{3} & {\cdot} & {\cdot} \\
\cline{1-3}
{\cdot} & {\cdot} & {\cdot} \\
\end{array}$}}$}
\ \ ,\ \ \ \ \ \  \ 
Q\ = \
\resizebox{.18\hsize}{!}{$
{\def\lr#1{\multicolumn{1}{|@{\hspace{.75ex}}c@{\hspace{.75ex}}|}{\raisebox{-.04ex}{$#1$}}}\raisebox{-.6ex}
{$\begin{array}{cccccc}
\cdot &\cdot & \cdot &\cdot &\cdot &\cdot \\
\cline{4-5}
\cdot & \cdot & \cdot & \lr{5} & \lr{6} \\
\cline{1-5}
\lr{4} & \lr{4} & \lr{4} \\
\cline{1-3}
\end{array}$}}$},$$
and  
$$ P \ = \ \sigma^{-6}(P')\ = \
\resizebox{.18\hsize}{!}{$
{\def\lr#1{\multicolumn{1}{|@{\hspace{.75ex}}c@{\hspace{.75ex}}|}{\raisebox{-.04ex}{$#1$}}}\raisebox{-.6ex}
{$\begin{array}{cccccc}
 &   & & & &  \\
\cline{4-6}
& & & \lr{\ov{3}} & \lr{\ov{3}} & \lr{\ov{2}} \\
\cline{1-6}
\lr{\ov{3}} & \lr{\ov{2}} & \lr{\ov{2}} & \lr{\ov{1}} & \lr{\ov{1}} & \lr{\ov{1}} \\
\cline{1-6}
\end{array}$}}$}\ .$$
Applying $\kappa^{-1}$ to the pair \vskip 2mm
$$ (P,Q) \ =  \left( \ \ 
\resizebox{.18\hsize}{!}{$
{\def\lr#1{\multicolumn{1}{|@{\hspace{.75ex}}c@{\hspace{.75ex}}|}{\raisebox{-.04ex}{$#1$}}}\raisebox{-.6ex}
{$\begin{array}{cccccc}
\cline{4-6}
& & & \lr{\ov{3}} & \lr{\ov{3}} & \lr{\ov{2}} \\
\cline{1-6}
\lr{\ov{3}} & \lr{\ov{2}} & \lr{\ov{2}} & \lr{\ov{1}} & \lr{\ov{1}} & \lr{\ov{1}} \\
\cline{1-6}
\end{array}$}}$}
\ \ \ ,\ \ \
\resizebox{.18\hsize}{!}{$
{\def\lr#1{\multicolumn{1}{|@{\hspace{.75ex}}c@{\hspace{.75ex}}|}{\raisebox{-.04ex}{$#1$}}}\raisebox{-.6ex}
{$\begin{array}{cccccc}
\cline{4-5}
 & & & \lr{5} & \lr{6} & \cdot \\
\cline{1-5}
\lr{4} & \lr{4} & \lr{4} & \cdot & \cdot & \cdot \\
\cline{1-3}
\end{array}$}}$}\ \ \right),
$$\vskip 2mm
\noindent where ${\rm sh}(P)=(6,3)^\pi$ and ${\rm sh}(Q)=\left((6,3)/(3,1)\right)^\pi$,
we have $(T,M)$ where \vskip 2mm
$$  
T \ = \ 
\resizebox{.09\hsize}{!}{$
{\def\lr#1{\multicolumn{1}{|@{\hspace{.75ex}}c@{\hspace{.75ex}}|}{\raisebox{-.04ex}{$#1$}}}\raisebox{-.6ex}
{$\begin{array}{ccc}
\cline{3-3}
& & \lr{\ov{3}} \\
\cline{1-3}
\lr{\ov{2}} & \lr{\ov{2}} & \lr{\ov{1}} \\
\cline{1-3}
\end{array}$}}$} 
\quad \in SST_{[\ov{3}]}((3,1)^\pi) , \quad \
M\ = \
\begin{bmatrix}
0 & 1 & 1 \\
1 & 0 & 0 \\
2 & 0 & 0 \\
\end{bmatrix}\ \in \M_{[\ov{3}]\times ([6]\setminus [3])},
$$ \vskip 2mm

\noindent with
$$ S^- = T^{\nw} =
\resizebox{.09\hsize}{!}{$
{\def\lr#1{\multicolumn{1}{|@{\hspace{.75ex}}c@{\hspace{.75ex}}|}{\raisebox{-.04ex}{$#1$}}}\raisebox{-.6ex}
{$\begin{array}{ccc}
\cline{1-3}
\lr{\ov{3}} & \lr{\ov{2}} & \lr{\ov{1}} \\
\cline{1-3}
 \lr{\ov{2}} & & \\
\cline{1-1}
\end{array}$}}$}\ .
$$\vskip 2mm

\noindent Therefore, we have a triple $(S^+,S^-,M)$ associated to $S$:
$$
(S^+,S^-,M)\ = \
\left(\ \ 
\resizebox{.09\hsize}{!}{$
{\def\lr#1{\multicolumn{1}{|@{\hspace{.75ex}}c@{\hspace{.75ex}}|}{\raisebox{-.04ex}{$#1$}}}\raisebox{-.6ex}
{$\begin{array}{ccc}
\cline{1-3}
\lr{5}& \lr{5} & \lr{6} \\
\cline{1-3}
\lr{6} & \lr{6}  \\
\cline{1-2} 
\end{array}$}}$}\ \ , \ \
\resizebox{.09\hsize}{!}{$
{\def\lr#1{\multicolumn{1}{|@{\hspace{.75ex}}c@{\hspace{.75ex}}|}{\raisebox{-.04ex}{$#1$}}}\raisebox{-.6ex}
{$\begin{array}{ccc}
\cline{1-3}
\lr{\ov{3}} & \lr{\ov{2}} & \lr{\ov{1}} \\
\cline{1-3}
 \lr{\ov{2}} & & \\
\cline{1-1}
\end{array}$}}$}\ \ , \ \
\begin{bmatrix}
0 & 1 & 1 \\
1 & 0 & 0 \\
2 & 0 & 0 \\
\end{bmatrix}\ \ 
\right).
$$\vskip 2mm

\noindent Finally, the corresponding ${\bf c}(S)=({\bf c}^J,{\bf c}_{J_1},{\bf c}_{J_2})\in \cB_{\Omega}$ in \eqref{eq:c} is given by \vskip 2mm
\begin{equation*}
\begin{split}
{\bf c}^J&=
\begin{bmatrix}
c_{34} & c_{35} & c_{36} \\
c_{24} & c_{25} & c_{26} \\
c_{14} & c_{15} & c_{16} \\
\end{bmatrix}
=\begin{bmatrix}
0 & 1 & 1 \\
1 & 0 & 0 \\
2 & 0 & 0 \\
\end{bmatrix}, 
\end{split}
\end{equation*}
and 
\begin{equation*}
\begin{split}
{\bf c}_{J_1}&= 
\begin{bmatrix}
c_{23} & c_{13} \\
       & c_{12}
\end{bmatrix}
=
\begin{bmatrix}
1 & 1 \\
  & 0
\end{bmatrix},\ \ \ \ \ 
{\bf c}_{J_2} = 
\begin{bmatrix}
c_{45} & c_{46} \\
       & c_{56}
\end{bmatrix}
=
\begin{bmatrix}
2 & 1 \\
  & 2
\end{bmatrix}.
\end{split}
\end{equation*}
}
\end{ex}\vskip 2mm

\begin{rem}{\rm 
Let $\Omega'$ be another quiver of type $A_{n-1}$ with a single sink. Using Theorem \ref{thm:main-1}, one can describe the transition map $R_{\Omega}^{\Omega'} : \cB_\Omega \rightarrow \cB_{\Omega'}$ as follows. 

Let ${\bf c}\in \cB_{\Omega}$ be given. There exist a pair of Young tableaux $(S^+,S^-)$ (but not necessarily unique) such that  
${\bf c}_{J_1}={\bf c}^-(S^-)$ and ${\bf c}_{J_2}={\bf c}^+(S^+)$. We can apply the inverse algorithm of \eqref{eq:embedding} to $(S^+,S^-,{\bf c}^J)$ to obtain $S\in SST_{[n]}(\la)$ for some $\la\in \cP_n$ such that each $\la_i-\la_{i+1}$ is sufficiently large. In fact, we obtain a unique {(marginally) large tableau} (see \cite{CL,HL}) corresponding to ${\bf c}$.  Let ${\bf c}'$ be the Lusztig datum of $S$ with respect to $\Omega'$, which is also obtained by the algorithm \eqref{eq:embedding}. Then we have ${\bf c}'=R_{\Omega}^{\Omega'} (\bf c)$. 

Note that if either one of $\Omega$ and $\Omega'$ is $\Omega^\pm$, then one may apply only Propositions \ref{prop:embedding Omega^+} and \ref{prop:embedding Omega^-} to have $R_{\Omega}^{\Omega'}$. We also refer the reader to \cite[Section 4]{BFZ} for a closed-form formula for $R_{\Omega}^{\Omega'}$, which is a tropicalization of a subtraction-free rational function connecting two parametrizations of a totally positive variety. It would be interesting to compare these two algortithms.}
\end{rem}


{\small
\begin{thebibliography}{HK}

\bibitem{BKT}
P. Baumann, J. Kamnitzer, P. Tingley, 
{\em Affine Mirkovi\'{c}-Vilonen polytopes}, Publ. Math. Inst. Hautes \'{E}tudes Sci. \textbf{120} (2014) 113--205. 
 
\bibitem{BFZ}
A. Berenstein, S. Fomin, A. Zelevinsky, {\em Parametrization of canonical bases and totally positive matrices}, Adv. in Math. \textbf{122} (1996) 49--149. 
 
 
\bibitem{CL} G. Cliff, {\em Crystal bases and Young tableaux}, J. Algebra \textbf{202} (1998) 10--35.
 
\bibitem{HK}
J. Hong, S.-J. Kang, {\em Introduction to Quantum Groups and Crystal Bases}, Graduate Studies in Mathematics 42,  Amer. Math. Soc., 2002.

\bibitem{HL} 
J. Hong, H.-M. Lee, {\em Young tableaux and crystal $B(\infty)$ for finite simple Lie
algebras}, J. Algebra \textbf{320} (2008) 3680--3693.

\bibitem{Ful}
W. Fulton, {\em Young tableaux, with Application to Representation
theory and Geometry}, Cambridge Univ. Press, 1997.



\bibitem{Kas91}
M. Kashiwara, {\em On crystal bases of the $q$-analogue of universal enveloping algebras}, Duke Math. J. \textbf{63} (1991) 465--516.

\bibitem{Kas93}
M. Kashiwara, {\em The crystal base and  Littelmann's refined Demazurecharacter formula}, Duke Math. J. \textbf{71} (1993) 839--858.

\bibitem{Kas95}
M. Kashiwara, \emph{On crystal bases}, Representations of groups,
CMS Conf. Proc., 16, Amer. Math. Soc., Providence, RI, (1995) 155--197.

\bibitem{KashNaka}
M. Kashiwara, T. Nakashima, {\em Crystal graphs for representations of the $q$-analogue of classical Lie algebras}, J. Algebra \textbf{165} (1994) 295--345.

\bibitem{Ki}
Y. Kimura, {\em Remarks on quantum unipotent subgroup and dual canonical basis}, preprint (2015) arXiv:1506.07912.

\bibitem{K09}
J.-H. Kwon,
{\em Demazure crystals of generalized Verma modules and a flagged RSK correspondence}, 
J. Algebra \textbf{322} (2009) 2150--2179.

\bibitem{K13}
J.-H. Kwon,
{\em RSK correspondence and classically irreducible Kirillov-Reshetikhin crystals}, 
J. Combin. Theory Ser. A \textbf{120} (2013) 433--452.

\bibitem{K14}
J.-H. Kwon, {\em Crystal bases of q-deformed Kac modules}, Int. Math. Res. Not. (2014)  512--550. 

\bibitem{Lu90}
G. Lusztig, {\em Canonical bases arising from quantized universal enveloping algebras},  J. Amer. Math. Soc.  \textbf{3}  (1990)  447--498.


\bibitem{Lu91}
G. Lusztig, {\em Quivers, perverse sheaves, and quantized enveloping algebras}, J. Amer. Math. Soc. \textbf{4} (1991) 365--421.

\bibitem{Lu93}
G. Lusztig, {\em Introduction to quantum groups}, Progress in Math. \textbf{110}, Birkh\"{a}user, 1993.

\bibitem{Lu96}
G. Lusztig, {\em Braid group action and canonical bases}, Adv. Math. \textbf{122} (1996) 237--261. 


\bibitem{Re}
M. Reineke,
{\em On the coloured graph structure of Lusztig's canonical basis}, Math. Ann. \textbf{307} (1997) 705--723.

\bibitem{S94}
Y. Saito, {\em  PBW basis of quantized universal enveloping algebras}, Publ. Res. Inst. Math. Sci. \textbf{30} (1994) 209--232.

\bibitem{S12}
Y. Saito, {\em Mirković-Vilonen polytopes and a quiver construction of crystal basis in type A}, Int. Math. Res. Not. IMRN (2012) 3877--3928.

\bibitem{SS}
B. E. Sagan, R. Stanley, {\em Robinson-Schensted algorithms for skew
tableaux}, J. Combin. Theory Ser. A \textbf{55} (1990) 161--193.

\bibitem{SST}
B. Salisbury, A. Schultze, P. Tingley, {\em Combinatorial descriptions of the crystal
structure on certain PBW bases}, preprint (2016) arXiv:1606.01978v2.

\bibitem{T}
T. Tanisaki, {\em Modules over quantized coordinate algebras and PBW-bases}, preprint (2014) arXiv:1409.7973.

\end{thebibliography}
}
\end{document}